\documentclass{amsart}
\usepackage{amssymb,pinlabel}

\def\co{\colon\thinspace}
\newcommand{\SL}{\widetilde{\mbox{\rm SL}}_2}
\newcommand{\EE}{\widetilde{\mbox{\rm E}}_2}
\newcommand{\SO}{\mbox{\rm SO}}
\newcommand{\SU}{\mbox{\rm SU}}
\newcommand{\oq}{\overline{q}}
\newcommand{\og}{\overline{g}}
\newcommand{\oz}{\overline{z}}
\newcommand{\ol}{\overline{\ell}}
\newcommand{\ou}{\overline{u}}
\newcommand{\tu}{\tilde{u}}
\newcommand{\otwo}{\overline{2}}
\newcommand{\obeta}{\overline{\beta}}
\newcommand{\odelta}{\overline{\delta}}
\newcommand{\opartial}{\overline{\partial}}
\newcommand{\bi}{\mathbf{i}}
\newcommand{\bj}{\mathbf{j}}
\newcommand{\bk}{\mathbf{k}}
\newcommand{\bv}{\mathbf{v}}
\newcommand{\bw}{\mathbf{w}}
\newcommand{\bz}{\mathbf{z}}
\newcommand{\bY}{\mathbf{Y}}
\newcommand{\bG}{\mathbf{G}}
\newcommand{\bZ}{\mathbf{Z}}
\newcommand{\R}{\mathbb{R}}
\newcommand{\C}{\mathbb{C}}
\newcommand{\D}{\mathbb{D}}
\newcommand{\E}{\mathbb{E}}
\newcommand{\HH}{\mathbb{H}}

\newcommand{\Z}{\mathbb{Z}}
\newcommand{\G}{\mathcal G}
\newcommand{\PSL}{\mbox{\rm PSL}}
\newcommand{\ip}{{\,\rule{2.3mm}{.2mm}\rule{.2mm}{2.3mm}\;\, }}

\newtheorem{theorem}{Theorem}
\newtheorem{lemma}[theorem]{Lemma}
\newtheorem{proposition}[theorem]{Proposition}
\newtheorem{corollary}[theorem]{Corollary}
\newtheorem{definition}[theorem]{Definition}
\newtheorem{remark}[theorem]{Remark}

\theoremstyle{definition}
\newtheorem*{remnot}{Remark on Notation}
\newtheorem*{acknowledgements}{Acknowledgements}

\begin{document}

\title{Contact Spheres and Hyperk\"ahler Geometry}

\author{Hansj\"org Geiges}
\address{Mathematisches Institut, Universit\"at zu K\"oln,
Weyertal 86--90, 50931 K\"oln, Germany.}
\email{geiges@math.uni-koeln.de}
\thanks{H. G. is partially supported by
         DFG grant GE 1245/1-2 within the framework of the
         Schwer\-punkt\-programm ``Globale Differentialgeometrie''.}

\author{Jes\'us Gonzalo P\'erez}
\address{Departamento de Matem\'aticas,
Universidad Aut\'onoma de Ma\-drid, 28049 Madrid, Spain.}
\email{jesus.gonzalo@uam.es}
\thanks{J. G. is partially supported by
         grants MTM2004-04794 and MTM2007-61982 from MEC Spain.}

\date{}

\begin{abstract}
\noindent A taut contact sphere on a $3$-manifold is a linear
$2$-sphere of contact forms, all defining the same volume form.
In the present paper we completely determine the moduli of taut
contact spheres on compact left-quotients of $\SU (2)$ (the only
closed manifolds admitting such structures). We also show that the
moduli space of taut contact spheres embeds into the moduli space of
taut contact circles.

This moduli problem leads to a new viewpoint on the Gibbons-Hawking
ansatz in hyperk\"ahler geometry. The classification of
taut contact spheres on closed $3$-manifolds includes
the known classification of
$3$-Sasakian $3$-manifolds, but the local
Riemannian geometry of contact spheres 
is much richer. We construct two examples
of taut contact spheres on open subsets of $\R^3$
with nontrivial local geometry; one from the Helmholtz equation on
the $2$-sphere, and one from the Gibbons-Hawking ansatz.
We address the Bernstein problem whether such examples can give rise
to complete metrics.
\end{abstract}

\maketitle



\section{Introduction}
We begin with the definition of our basic objects of interest.
Recall that a {\bf contact form} on a $3$-manifold is a
differential $1$-form $\alpha$ such that $\alpha\wedge d\alpha
\neq 0$.

\begin{definition}
A {\bf contact sphere} is a triple of $1$-forms
$(\alpha_1, \alpha_2, \alpha_3)$ on a $3$-manifold
such that any non-trivial linear
combination of these forms (with constant coefficients) is a
contact form.
\end{definition}

In other words, we require that the $3$-form
\[ (\lambda_1\alpha_1+\lambda_2\alpha_2+\lambda_3\alpha_3)
  \wedge (\lambda_1\, d\alpha_1+\lambda_2\, d\alpha_2+\lambda_3\, d\alpha_3) \]
be nowhere zero, i.e.\ a volume form, for any $\lambda_1,\lambda_2,\lambda_3
\in \R$ with $\lambda_1^2+\lambda_2^2+\lambda_3^2\neq 0$.
The name `contact sphere' derives from the fact that it suffices to
check this condition for points $(\lambda_1,\lambda_2,\lambda_3)$
on the unit sphere $S^2\subset \R^3$.

\begin{definition}
A contact sphere $(\alpha_1,\alpha_2,\alpha_3)$ is called {\bf taut}
if the contact form $\lambda_1\alpha_1+\lambda_2\alpha_2+\lambda_3
\alpha_3$ defines the same volume form for all $(\lambda_1,\lambda_2,
\lambda_3)\in S^2$.
\end{definition}

The requirement for a contact sphere to be taut is equivalent to
the system of equations (for $i\neq j$)
\begin{eqnarray*}
\alpha_i\wedge d\alpha_i & = & \alpha_j\wedge d\alpha_j\neq 0,\\
\alpha_i\wedge d\alpha_j & = & -\alpha_j\wedge d\alpha_i.
\end{eqnarray*}
A straightforward calculation shows that one can then find a $1$-form
$\beta$ and a nowhere zero function $\Lambda$ such that
\begin{equation}
\label{eqn:structure}
d\alpha_i=\beta\wedge\alpha_i +\Lambda\,\alpha_j\wedge\alpha_k,
\end{equation}
where $(i,j,k)$ runs over the cyclic permutations of $(1,2,3)$.
Notice that $\Lambda$ is defined by $\alpha_i\wedge d\alpha_i=\Lambda\,
\alpha_1\wedge \alpha_2\wedge\alpha_3$.

The analogous structure of a (taut) contact circle, defined in terms of
two contact forms $(\alpha_1,\alpha_2)$, was studied in our
previous papers \cite{gego95}, \cite{gego97}, \cite{gego02}. In \cite{gego95}
we gave a complete classification of the closed, orientable
$3$-manifolds\footnote{Our manifolds are always understood to be
connected.} that
admit a taut contact circle or a taut contact sphere:

\begin{theorem}
\label{thm:classold}
Let $M$ be a closed $3$-manifold.
\begin{itemize}
\item[{\rm (a)}] $M$ admits
a taut contact circle
if and only if $M$ is diffeomorphic to a quotient of the Lie group
$\G$ under a discrete subgroup $\Gamma$ acting by left multiplication,
where $\G$ is one of the following.
\begin{itemize}
\item[{\rm (i)}] $S^{3}=\SU (2)$, the universal cover of $\SO(3)$.

\item[{\rm (ii)}] $\SL$, the universal cover of $\PSL_{2}\R$.

\item[{\rm (iii)}] $\EE$, the universal cover of the Euclidean
group (that is, orientation preserving isometries of~$\R^2$).
\end{itemize}

\item[{\rm (b)}] $M$ admits a taut contact sphere if and only if it is
diffeomorphic to a left-quotient of $\SU (2)$.
\end{itemize}
\end{theorem}

In the course of this  paper we shall present a new proof,
more self-contained than the one given in~\cite{gego95},
of the fact that the universal cover of a closed $3$-manifold
admitting a taut contact sphere is diffeomorphic to~$S^3$.

In \cite{gego97} we showed that every closed, orientable $3$-manifold
admits a (non-taut) contact circle, and we gave examples of contact spheres.
For instance, $S^1\times S^2\subset S^1\times \R^3$, described in terms
of coordinates $(\theta ,x,y,z)$, does not admit any taut contact
circles by Theorem~\ref{thm:classold}, but it admits the contact
sphere
\begin{eqnarray*}
\alpha_1 & = & x\, d\theta +y\, dz-z\, dy,    \\
\alpha_2 & = & y\, d\theta +z\, dx -x\, dz,   \\
\alpha_3 & = & z\, d\theta +x\, dy -y\, dx.
\end{eqnarray*}

In \cite{gego02} we described deformation spaces for taut contact circles
and gave a complete classification of taut contact circles. The present
paper achieves the corresponding classification for taut contact spheres.

The investigation of taut contact spheres amounts to a systematic
study of hyperk\"ahler metrics with a homothety,
in a sense made precise below.
Constructions of complete hyperk\"ahler metrics with
translational {\em invariance} have played a prominent role in general
relativity and supersymmetric field theories, beginning with the
Gibbons-Hawking ansatz~\cite{giha78}. This ansatz will be discussed
in Section~\ref{section:GH}
in the context of explicit constructions of taut contact spheres
with nontrivial local geometry. See~\cite{basf97} for a
fairly recent discussion of several constructions related
to the Gibbons-Hawking ansatz
(Eguchi-Hanson metric, Taub-NUT metric, Atiyah-Hitchin metric).

Hyperk\"ahler metrics with a {\em homothety} are
equally important
in physical applications, see~\cite{ctv96} and~\cite{giry98}. The
latter pays particular attention to homotheties that are
hypersurface orthogonal (such homotheties are called
{\it dilatations} in~\cite{giry98}). This is equivalent to saying that
the hyperk\"ahler metric on $U\times \R$, where $U$ is some
$3$-manifold, is the cone metric
over a $3$-Sasakian metric on $U$, cf.~\cite{bar93}. For more
general information on Sasakian and in particular $3$-Sasakian
geometry see the definitive survey~\cite{boga99} or the
monograph~\cite{boga08}; some of the
definitions will be recalled in Section~\ref{section:Sasakian}.
There we use our methods to recover the classification
of the closed $3$-manifolds admitting $3$-Sasakian structures.
In contrast with $3$-Sasakian structures,
taut contact spheres do not, in general, give rise to a cone
metric on $U\times \R$ (in other words, the relevant
homothety is not a dilatation; such general homotheties also
appear in~\cite{ctv96}). This implies that taut contact spheres
are definitely more general than $3$-Sasakian structures. We
elaborate on this point in Section~\ref{section:Sasakian}.

By comparison, K\"ahler metrics on $U\times \R$ admitting a
dilatation correspond to a Sasakian structure on~$U$. (For a classification
of the closed $3$-dimensional manifolds admitting Sasakian structures
see~\cite{geig97}.) In~\cite{gkn00} it is shown that if the Sasakian
analogue of the K\"ahler potential satisfies a Monge-Amp\`ere equation,
then the metric on $U$ is Sasakian-Einstein. This happens in particular if
$U\times \R$ is Ricci-flat, cf.~\cite{giry98}.

A taut contact sphere on a closed manifold $M$ always gives rise
to a {\em flat} hyperk\"ahler metric on $M\times\R$
(Theorem~\ref{thm:flat}). This must be read as a
global rigidity phenomenon, because the theorem fails
for open $3$-manifolds. In Section~\ref{section:Helmholtz}
we use a Monge-Amp\`ere equation to construct a
taut contact sphere on an open subset $U$ of $\R^3$
giving rise to a non-flat hyperk\"ahler metric on $U\times\R$.
In the appendix we
use a contact transformation to relate this construction
to the Helmholtz equation on the $2$-sphere. In Section~\ref{section:GH}
we use the Gibbons-Hawking ansatz to give an even simpler
construction of a non-flat example. In Section~\ref{section:Bernstein}
we discuss the question, known as a {\em Bernstein problem}, whether
such non-flat examples can give rise to complete metrics. The answer to
this question depends on
the choice of one of the two natural metrics associated with a taut
contact sphere (see Definition~\ref{defn:short-long},
Theorem~\ref{thm:Bernstein}, and the comments following it).

This paper supersedes our preprint
``Contact spheres and quaternionic structures''.
\section{Statement of results}
We now describe in outline some of the main results of the present paper.
Our notational convention throughout will be that $M$ denotes a
closed, orientable $3$-manifold; $U$ will denote a $3$-manifold
(without boundary) that need not be compact.

The relation
\[ (\alpha_1,\alpha_2,\alpha_3)\sim (v\alpha_1,v\alpha_2,v\alpha_3)\;\;
\mbox{\rm for some smooth function}\;\; v\co U\rightarrow \R^+\]
is easily seen to be an equivalence relation within the set of
(taut) contact spheres.

\begin{definition}
Two (taut) contact spheres are {\bf conformally equivalent} if one is
obtained from the other by multiplying each contact form by the same
positive function.
\end{definition}

\begin{definition}
We call a contact sphere {\bf naturally ordered} if
$\alpha_i\wedge d\alpha_i$ is a positive multiple of
$\alpha_1\wedge\alpha_2\wedge\alpha_3$.
\end{definition}

Throughout this paper we shall assume our contact spheres to satisfy this
condition.

Since $\alpha_i\wedge d\alpha_i$ and $\alpha_1\wedge\alpha_2\wedge
\alpha_3$ scale with the second and third power of~$v$,
respectively, it is obvious that every conformal equivalence class
of naturally ordered 
taut contact spheres contains, for any $c\in\R^+$,
a unique representative satisfying
\begin{equation}
\label{eqn:normalised}
\alpha_i\wedge d\alpha_i= c\,\alpha_1\wedge\alpha_2\wedge\alpha_3,\;\;
i=1,2,3.
\end{equation}

\begin{definition}
We call a taut contact sphere $(\alpha_1,\alpha_2,\alpha_3)$
{\bf $c$-normalised} if it satisfies
equation~(\ref{eqn:normalised}).
\end{definition}

\begin{remark}
{\rm
This condition is equivalent
to $\Lambda\equiv c$ in equation~(\ref{eqn:structure}); beware that
$\beta$ in that equation is not an invariant of the conformal
equivalence class.
}
\end{remark}

It is implicit in \cite{geig96} and follows by a simple extension of the
ideas from \cite{gego95} that a taut contact sphere on $U$ gives rise
to a hyperk\"ahler structure on $U\times \R$. In
Section~\ref{section:hyper} we analyse this situation a little more
carefully. One of the results proved there is the following.

\begin{proposition}
\label{prop:hyper1}
A contact sphere on $U$ determines an oriented conformal structure
on $U\times \R$. A naturally ordered taut contact sphere on $U$ determines
a hyperk\"ahler structure on $U\times \R$. Conformally
equivalent (taut) contact spheres determine isomorphic conformal
(resp.\ hyperk\"ahler) structures.
\end{proposition}

As we shall see in Section~\ref{section:hyper}, the
hyperk\"ahler structure $(g,J_1,J_2,J_3)$ on $U\times
\R$ induced by a taut contact sphere $(\alpha_1,\alpha_2,\alpha_3)$
on $U$ is given by the equations
\[ -g(\cdot ,J_i\cdot )=d(e^t\alpha_i)=:\Omega_i,\; i=1,2,3,\]
where $t$ denotes the $\R$-coordinate, and the complex structures $J_i$
are $\partial_t$-invariant. Often we write
the hyperk\"ahler structure as the triple $(\Omega_1,\Omega_2,\Omega_3)$
of symplectic forms.

These symplectic forms are homogeneous of degree $1$ with respect to
the vector field $\partial_t$, that is,
$L_{\partial_t}\Omega_i=\Omega_i$. The hyperk\"ahler metric $g$
has the same homogeneity.

For reversing this construction, it is useful to
make the following definition.

\begin{definition}
A vector field $Y$ on a hyperk\"ahler manifold is called
{\bf tri-Liouville} if it is a Liouville vector field
for every parallel self-dual $2$-form~$\Omega$, that is,
$L_X\Omega =\Omega$.
\end{definition}

Notice that a tri-Liouville vector field is automatically
homothetic for the hyperk\"ahler metric, but the converse is
not true.

The construction of taut contact spheres in this paper uses two
ingredients: a hyperk\"ahler metric $g$
and a nowhere zero tri-Liouville vector field~$Y$. For any
conformal basis $(\Omega_1,\Omega_2,\Omega_3)$ of parallel
self-dual $2$-forms, the corresponding taut contact sphere
is defined on any transversal to the flow of $Y$ by restricting the
triple of $1$-forms $(i_Y\Omega_1,i_Y\Omega_2,i_Y\Omega_3)$ to
that transversal. For the formal statement see
Proposition~\ref{prop:hyper2}.

We also show that for a naturally ordered
taut contact sphere $(\alpha_1,\alpha_2,
\alpha_3)$ on $U$ there is the following {\em pointwise} model
on $U\times \R$ for the triple of symplectic forms
$(\Omega_1,\Omega_2,\Omega_3)$,
expressed in quaternionic notation: At any point $x\in U\times \R$,
there is a quaternionic coordinate $dq_x$ for the tangent space
$T_x(U\times \R)$ such that
\[ d(e^t(\bi\alpha_1+\bj\alpha_2+\bk\alpha_3))_x=- dq_x
\wedge d\oq_x.\]

The key to the classification of taut contact spheres is then the
following statement.

\begin{theorem}[Global rigidity]
\label{thm:flat}
The hyperk\"ahler metric on $M\times \R$ induced by a taut
contact sphere on a closed $3$-manifold $M$ is flat.
\end{theorem}

\begin{proof}
Hyperk\"ahler metrics are always Ricci flat~\cite[14.13]{bess87},
and any Ricci flat K\"ahler manifold of complex dimension $2$
is anti-self-dual~\cite{ahs78}, that is, the self-dual
part $W^+$ of the Weyl tensor of the metric $g$ vanishes.
Since the Weyl tensor is an invariant of the conformal class of a
metric~\cite[1.159]{bess87}, $W^+$ also vanishes for the
metric $g/g(\partial_t,\partial_t)$,
which descends to the quotient $M\times S^1$ of
$M\times \R$ under the map $(p,t)\mapsto (p,t+1)$, say.
(In fact, the hyperhermitian structure $(g/g(\partial_t,\partial_t),
J_1,J_2,J_3)$
descends to that quotient, and one may also appeal to a result
of Boyer~\cite{boye88} saying that a hyperhermitian metric on
a $4$-manifold is anti-self-dual.)

Then the signature formula for $(M\times S^1,
g/g(\partial_t,\partial_t))$ yields
\begin{eqnarray*}
\tau (M\times S^1) & = & \frac{1}{12\pi^2}\int_{M\times S^1}
                          \| W^+\|^2-\|W^-\|^2\\
                   & = & -\frac{1}{12\pi^2}\int_{M\times S^1}
                           \|W^-\|^2.
\end{eqnarray*}
But $\tau(M\times S^1)=0$ for purely topological reasons. So the
Weyl tensor $W=W^++W^-$ vanishes for~$g/g(\partial_t,\partial_t)$.
Again appealing to the conformal invariance of the Weyl tensor, we deduce
that it also vanishes for~$g$.
For Ricci flat metrics (thus, in particular, for the hyperk\"ahler
metric~$g$) the full curvature tensor equals its
Weyl part~\cite[1.116]{bess87}. This proves the theorem.
\end{proof}

This theorem implies that on a closed manifold
the {\em pointwise} model above for the
triple of symplectic forms $d(e^t\alpha_i)$, $i=1,2,3$,
coming from a taut contact sphere, is actually a {\em local}
model, since the three symplectic forms and the hypercomplex
structure $(J_1,J_2,J_3)$ are all parallel with respect to the
flat hyperk\"ahler metric~$g$.

It is then not very difficult, using properties of the $\partial_t$-flow,
to derive the following classification statement. The proof will be
given in Section~\ref{section:class}.

\begin{theorem}
\label{thm:class}
If $M$ is diffeomorphic to a lens space $L(m,m-1)$, including
the $3$-sphere $L(1,0)$,
then the (naturally ordered) taut contact spheres on $M$,
up to diffeomorphism
and conformal equivalence, are given by the following family of
$\Z_m$-invariant quaternionic $1$-forms on $S^3\subset \HH$,
\[ \bi\alpha_1+\bj\alpha_2+\bk\alpha_3=\frac{1}{2}(dq\cdot\oq -
q\cdot d\oq )-\nu\, d(q\bi\oq ),\;\;\nu\in\R ,\]
where $\nu$ and $-\nu$ determine equivalent structures.

If $M$ is diffeomorphic to $\Gamma\backslash \SU (2)$ with $\Gamma
\subset \SU (2)$ a non-abelian group, there is a unique equivalence
class of taut contact
spheres on~$M$, described by the formula above with~$\nu =0$.

All these taut contact spheres are homogeneous under a natural
$\SO (3)$-action. In particular, all great circles in a given
taut contact sphere are isomorphic taut contact circles.
\end{theorem}

Note that the manifolds listed in this theorem exhaust all the possible
left-quotients of $\SU (2)$, cf.~\cite{gego95}. The $\Z_m$-action
on $S^3$ that produces the quotient $L(m,m-1)$ is generated by
right multiplication with $\cos(2\pi /m)+\bi\sin (2\pi /m)$.

\begin{remark}
{\rm
As we shall see in the construction of the moduli space of taut contact
spheres, all taut contact spheres on a given closed $3$-dimensional
manifold yield the same $4$-dimensional metric $g$ up to global
isometry. So it is not this induced metric which determines the
modulus, but in fact the different possible homothetic vector
fields~$\partial_t$.
In the cases where non-trivial moduli exist,
neither the vector field $\partial_t$ nor the
function $e^t$ (which turns out to equal $g(\partial_t,\partial_t)$)
are unique for that~$g$.
}
\end{remark}

With our coordinate
conventions in Section~\ref{section:hyper}
below, which seem natural in that context, {\em left} multiplication on
$\C^2$ by
\[ \left( \begin{array}{cc} a&\; -\overline{b}\\ b&\;\overline{a}\end{array}
\right) \in\SU (2),\]
where $a=a_1+\bi a_2$, $b=b_1+\bi b_2$, $|a|^2+|b|^2=1$, corresponds to
{\em right} multiplication on $\HH$ by the unit quaternion
$u=a_1+\bi a_2+\bj b_1+\bk b_2$.

The quaternionic $1$-form $dq\cdot \oq -q\cdot d\oq$ is
invariant under this right multiplication $q\mapsto qu$ by unit
quaternions~$u$, and therefore
descends to all left-quotients of~$\SU (2)$. The $1$-form
$d(q\bi\oq )$ is invariant under right multiplication by
unit quaternions of the special form $a_1+\bi a_2$, hence it
descends to all {\em abelian} quotients.

Here are some details about what we mean by `homogeneity under
a natural $\SO (3)$-action' in Theorem~\ref{thm:class}.
The action of $\SO (3)$ on $\bi\alpha_1+\bj\alpha_2+\bk\alpha_3$
rotates the contact forms and
is given by conjugating the quaternionic $1$-form by elements
$u\in S^3\subset\HH$.
This action is induced from the $S^3$-action $q\mapsto uq$ on~$\HH$,
and that latter action cannot be replaced by an $\SO (3)$-action. This amounts
to a spinor phenomenon, see also the proof of Proposition~\ref{prop:planes}.
Since this left multiplication by unit quaternions commutes
with all quaternionic right multiplications, it descends to
all left-quotients of~$\SU (2)$.

\begin{remark}
{\rm
With $q=x_0+\bi x_1+\bj x_2+\bk x_3$ we obtain the following real
expression for the taut contact sphere in Theorem~\ref{thm:class}
corresponding to $\nu =0$, where $(i,j,k)$ runs over the cyclic permutations
of $(1,2,3)$:
\[ \alpha_i=x_0\, dx_i-x_i\, dx_0+x_j\, dx_k-x_k\, dx_j.\]
Notice that it 
satisfies (\ref{eqn:structure}) with $\beta\equiv 0$ and $\Lambda\equiv 2$.
In particular, it is a $2$-normalised taut contact sphere.

Restricting this to the hyperplane $\{ x_0=1\}$, we obtain a
simple expression for a taut contact sphere on~$\R^3$:
\[ \alpha_i=dx_i+x_j\, dx_k-x_k\, dx_j. \]
See Proposition~\ref{prop:hyper2} for the principle behind this observation.
}
\end{remark}

One might suspect that taut contact spheres constitute such a rigid
structure that Theorem~\ref{thm:flat} would also hold locally and
for open manifolds. However, this turns out to be false, even conformally.

\begin{theorem}
\label{thm:nonflat}
There are examples of taut contact spheres on open domains $U$ in
$\R^3$ inducing a metric on $U\times \R$ that is not
conformally flat.
\end{theorem}

This theorem will be proved in Section~\ref{section:nonflat}, where
we present two methods for constructing such examples. Our
first construction, in Section~\ref{section:Helmholtz},
starts from the observation that
(locally) the conditions for a taut contact sphere lead
to a complex Monge-Amp\`ere equation.
After imposing an additional homogeneity amounting to the existence
of a tri-Hamiltonian vector field, we are able to find concrete
solutions of that equation whose associated $4$-dimensional metric is
non-flat. In fact, subject to this extra homogeneity, the
Monge-Amp\`ere equation can be linearised to yield
the Helmholtz equation $\Delta u+2u=0$ for the Laplacian
$\Delta$ on the $2$-sphere. Solutions of that Helmholtz equation
then give rise to taut contact spheres with a tri-Hamiltonian
symmetry.

Our second construction, in Section~\ref{section:GH}, starts from
the Gibbons-Hawking ansatz, which is essentially a construction
of an $\R$-invariant hyperk\"ahler metric on $U_0\times\R$,
starting from a harmonic function on an open subset $U_0$ of
$\R^3$. For an appropriate choice
of such a harmonic function, one obtains a non-flat hyperk\"ahler metric
giving rise to a taut contact sphere on a suitable hypersurface
$U\subset U_0\times\R$. Beware that the $\R$-factor in this splitting
$U_0\times\R$ is not the one corresponding to the tri-Liouville vector
field.

\vspace{2mm}

Taut contact spheres come associated with two natural metrics on
the $3$-manifold.
In order to describe these, we observe that the
$2$-forms $\Omega_i=d(e^t\alpha_i)$ can be written with the help
of the structure equation~(\ref{eqn:structure}) as follows:
\[ \Omega_i=e^t\bigl( \Lambda^{-1/2}(dt+\beta )\wedge \Lambda^{1/2}\alpha_i+
\Lambda^{1/2}\alpha_j\wedge \Lambda^{1/2}\alpha_k\bigl) .\]
So the hyperk\"ahler metric $g$ is given by
\begin{equation}
\label{eqn:metric}
g=e^t\bigl( \Lambda^{-1}(dt+\beta )^2+\Lambda\,
(\alpha_1^2+\alpha_2^2+\alpha_3^2)\bigr) .
\end{equation}
In particular, if the contact sphere is $1$-normalised, we have
\begin{equation}
\label{eqn:metric-normal}
g=e^t\bigl( (dt+\beta )^2+\alpha_1^2+\alpha_2^2+\alpha_3^2\bigr) .
\end{equation}
This motivates the following definition.

\begin{definition}
\label{defn:short-long}
The {\bf short metric} associated with a $1$-normalised taut contact sphere
$(\alpha_1,\alpha_2,\alpha_3)$ on $U$ is the metric
\[ g_s=\alpha_1^2+\alpha_2^2+\alpha_3^2.\]
The {\bf long metric} associated with this contact sphere is
\[ g_l=\beta^2+\alpha_1^2+\alpha_2^2+\alpha_3^2.\]
\end{definition}

Observe that $g_l$ is simply the restriction of $g$ to
$U\equiv U\times\{ 0\}$, so from the viewpoint of hyperk\"ahler
geometry this is the more natural metric to consider.

In either of the above-mentioned constructions of taut contact spheres
on $U$ giving rise to a non-flat
hyperk\"ahler metric on $U\times\R$, the induced
long metric $g_l$ on $U$ is incomplete, and so is, {\em a fortiori},
the short metric~$g_s$. In Section~\ref{section:Bernstein}
we raise the question whether one can find such examples where~$g_l$,
at least, is complete. This type
of question is known as a Bernstein problem~\cite{cala70}.
Concerning~$g_s$, we provide a partial
answer to this problem for contact spheres with
additional symmetries. For
$g_l$ there is a positive answer to the Bernstein problem, even subject
to the additional symmetry requirement:

\begin{theorem}
\label{thm:Bernstein}
Let $(\alpha_1,\alpha_2,\alpha_3)$ be a $1$-normalised
taut contact sphere on $U$ giving
rise to the short metric $g_s$ and the long metric $g_l$ on~$U$.

(a) If $g_s$ is complete and
admits a non-trivial Killing vector field that preserves
$(\alpha_1,\alpha_2,\alpha_3)$, then $U$ is necessarily compact, and hence
$(\alpha_1,\alpha_2,\alpha_3)$
belongs to the family described in Theorem~\ref{thm:class}
(and in particular gives rise to a flat hyperk\"ahler metric).

(b) There are examples of $S^1$-invariant taut contact spheres
on $U=\D\times S^1$, where $\D$ is the open unit disc,
for which $g_l$ is complete. The induced hyperk\"ahler metrics
on $U\times\R$ are not flat.
\end{theorem}

Part (a) will be proved in Section~\ref{section:Bernstein}.
Theorem~\ref{thm:Bernstein-Killing} in that section is
a more explicit reformulation of this part.

We reserve the proof of part (b) for a forthcoming paper.
It turns out that in the $\R$- or $S^1$-invariant context one can use the
Gibbons-Hawking ansatz in order to develop a complete theory
of such contact circles. An infinite-dimensional family
of examples giving rise to
complete long metrics is then found
with the help of Blaschke products on $\D\subset\C$.
\section{Hyperk\"ahler linear algebra}
\label{section:hyper}
In this section we discuss the linear algebraic aspects of contact
spheres, leading to a proof of Proposition~\ref{prop:hyper1}.
This prepares the ground for the proof
of Theorem~\ref{thm:class}.

Let $(\alpha_1,\alpha_2,\alpha_3)$ be a contact sphere on a
$3$-manifold~$U$. This gives rise to the symplectic forms
$\Omega_i=d(e^t\alpha_i)$, $i=1,2,3$, on $U\times \R$. At
any point $x$ of $U\times \R$, these symplectic forms
span a definite $3$-plane in the space $\bigwedge^2T^*_x(U\times
\R)$ of skew-symmetric bilinear forms on the tangent
space $T_x(U\times \R)$. If the contact sphere is taut,
we have in addition the identities (for $i\neq j$)
\begin{eqnarray*}
\Omega_i\wedge\Omega_i & = & \Omega_j\wedge\Omega_j\neq 0,\\
\Omega_i\wedge\Omega_j & = & 0.
\end{eqnarray*}

First we are going to study the linear algebra of this situation. Thus,
let $V_4$ be a $4$-dimensional real vector space and write $V_6=
\bigwedge^2V_4^*$. Consider the quadratic form
\[ Q\co V_6\longrightarrow \R ,\;\; Q(A)=A\wedge A,\]
of signature $(3,3)$. We call a triple $(A_1,A_2,A_3)$ of elements
of $V_6$ a {\bf symplectic triple} on $V_4$ if it spans a definite $3$-plane
for $Q$ in $V_6$, and a {\bf conformal symplectic triple} if the stronger
condition
\begin{eqnarray*}
A_i\wedge A_i & = & A_j\wedge A_j\neq 0,\\ A_i\wedge A_j & = & 0,
\end{eqnarray*}
is satisfied for $i\neq j$, cf.~\cite{geig96}. The same terminology will be used
for triples of symplectic forms $(\Omega_1,\Omega_2,\Omega_3)$ on a
$4$-manifold as described above.

\begin{remnot}
In the sequel, any equation (or other statement) involving the indices
$i,j,k$ is meant to be read as three equations, with $(i,j,k)$ ranging over
the cyclic permutations of $(1,2,3)$. We write bold face $\bi$ for
$\sqrt{-1}\in\C$, and bold face $\bi ,\bj ,\bk$ for the standard
quaternionic units with $\bi\bj =\bk$. The relation between real, complex, and
quaternionic coordinates will be given by
\[ z_1=x_0+\bi x_1,\;\; z_2=x_2+\bi x_3;\]
\[ q=x_0+\bi x_1+\bj x_2+\bk x_3 =z_1+z_2\bj.\]
\end{remnot}

For $J$ a complex structure on $V_4$, we denote by $\bigwedge ^{(p,q)}_J$
the space of exterior forms of type $(p,q)$ on~$V_4$. As shown
in~\cite{gego95}, a conformal symplectic couple $(A_i,A_j)$ determines
a unique complex structure $J_k$ on $V_4$ for which $A_i+\bi A_j\in
\bigwedge^{(2,0)}_{J_k}$. Thus, a conformal symplectic triple
$(A_1,A_2,A_3)$ induces three complex structures $J_1,J_2,J_3$.
In~\cite{geig96} it was shown that
these complex structures satisfy the quaternionic identities as well as
the relations
\[ g(\bv ,\bw )=A_i(\bv ,J_i\bw ), \;\; i=1,2,3, \;\;\mbox{\rm for all}\;
\bv ,\bw\in V_4,\]
for a unique definite symmetric bilinear form~$g$.
We call a conformal symplectic triple {\bf naturally ordered} if
this $g$ is {\em positive} definite. Notice that the sign of $g$
is well defined, since the $(-J_1,-J_2,-J_3)$ do not satisfy the
quaternionic identities. (In particular, hyperk\"ahler structures
are always naturally ordered.)

An alternative way to define this definite bilinear form $g$ is
via the identity
\begin{equation}
\label{eqn:cubic}
(\bv\ip A_1)\wedge (\bv\ip A_2)\wedge (\bv\ip A_3)=
\frac{1}{2}\, g(\bv ,\bv )\,  \bv\ip(A_i^2)\;\;\mbox{\rm for all}\;
\bv\in V_4.
\end{equation}
This is an obvious consequence of the $\SO (3)$-homogeneous
normal form for conformal
symplectic triples discussed in the next two propositions.

\begin{proposition}
\label{prop:pointmodel}
Let $(A_1,A_2,A_3)$ be
a naturally ordered conformal symplectic triple
on~$V_4$. Then there are real linear coordinates $dx_0,dx_1,dx_2,dx_3$
on $V_4$ such that
\[ A_i=dx_0\wedge dx_i+dx_j\wedge dx_k.\]
In terms of the corresponding complex and quaternionic coordinates we have
\begin{eqnarray*}
A_1 & = & \frac{\bi}{2}(dz_1\wedge d\oz_1+dz_2\wedge d\oz_2),\\
A_2+\bi A_3 & = & dz_1\wedge dz_2,
\end{eqnarray*}
and
\[ \bi A_1+\bj A_2+\bk A_3 = -\frac{1}{2}dq\wedge d\oq .\]
\end{proposition}

\begin{remark}
{\rm
Because of the non-commutativity of $\HH$ some care is necessary
in interpreting the wedge product of $\HH$-valued $1$-forms
$\alpha ,\beta$ on a vector space~$V$. Our convention is to read it as
\[ (\alpha\wedge\beta )(\bv ,\bw )=\alpha (\bv )\beta (\bw )-
\alpha (\bw )\beta (\bv )
\;\;\mbox{\rm for all}\; \bv ,\bw\in V.\] This ensures $\alpha \wedge
q\beta=\alpha q\wedge\beta$, and that
$\alpha\wedge\overline{\alpha}$ is always purely imaginary.
}
\end{remark}

\begin{proof}[Proof of Proposition~\ref{prop:pointmodel}]
From $0\neq A_2+\bi A_3\in\bigwedge^{(2,0)}_{J_1}$ and $A_1\wedge
(A_2+\bi A_3)=0$ we conclude that the $(0,2)$-part of $A_1$ with
respect to $J_1$ is zero. The form $A_1$ being real, its
$(2,0)$-part must also vanish, hence $A_1\in
\bigwedge^{(1,1)}_{J_1}$.

Let $\{\ell_1',\ell_2'\}$ be a basis for $\bigwedge^{(1,0)}_{J_1}$ and
set $\ell '=\left( \begin{array}{c} \ell_1'\\ \ell_2'\end{array} \right)$.
Then we can write
\[ A_1= \bi (\ell ')^T\wedge {\mathbf A}_1 \overline{\ell}',\]
with ${\mathbf A}_1$ a hermitian $(2\times 2)$-matrix. By our assumption
on $(A_1,A_2,A_3)$ being naturally ordered, the matrix ${\mathbf A}_1$
is positive definite. Hence there is a matrix ${\mathbf C}\in
\mbox{\rm GL}_2(\C )$ such that
\[ {\mathbf C}^T{\mathbf A}_1{\mathbf C}=\left( \begin{array}{cc}
1/2 &\; 0\\ 0 &\; 1/2 \end{array} \right ) ,\]
and so in terms of the basis $\{ \ell_1,\ell_2\}$
for $\bigwedge^{(1,0)}_{J_1}$
defined by $\ell '={\mathbf C}\ell$ we have
\[ A_1=\frac{\bi}{2}(\ell_1\wedge\ol_1+\ell_2\wedge\ol_2)\]
and
\[ A_2+\bi A_3= c\,\ell_1\wedge \ell_2\]
for some $c\in \C$. We then find
\begin{eqnarray*}
|c|^2\ell_1\wedge \ell_2\wedge\ol_1\wedge\ol_2 & = & (A_2+\bi A_3)
                                    \wedge (A_2-\bi A_3)\\
    & = & A_2^2+A_3^2\, =\,  2A_1^2\\
    & = & -\ell_1\wedge\ol_1\wedge \ell_2\wedge\ol_2,
\end{eqnarray*}
from which we conclude $|c|=1$. The linear complex coordinates
$z_1,z_2$ corresponding to $\{ c\,\ell_1,\ell_2\}$ then give the
desired complex normal form. The real and quaternionic normal
forms can be derived easily from the complex one.
\end{proof}

\begin{remark}
{\rm
\label{rem:metric}
(1) Notice that in terms of these pointwise
coordinates we have $g=dx_1^2+dx_2^2+dx_3^2+dx_4^2$ and $(J_1,J_2,J_3)
=(\bi ,\bj ,\bk )$. Moreover, we recognise the $3$-plane in $V_6$ spanned by
$A_1,A_2,A_3$ as the space of self-dual $2$-forms for the metric $g$ and
the orientation of $V_4$ defined by $A_i\wedge A_i$. The length of
the $A_i$ equals~$\sqrt{2}$.

(2) Here is another characterisation of taut contact spheres that
can be read off from the preceding proposition: The purely
imaginary $1$-form $\alpha =\bi\alpha_1+\bj\alpha_2+\bk\alpha_3$
defines a taut contact sphere if and only if $\alpha_1$ is a
contact form and at each point $x$ of the manifold there is an
$\HH$-valued linear form $\beta_x$ on the tangent space
at~$x$ such that $d(e^t\alpha
)_x=\pm\beta_x\wedge\overline{\beta}_x$. The case $d(e^t\alpha
)_x=-\beta_x\wedge\overline{\beta}_x$ corresponds to
$(\alpha_1,\alpha_2,\alpha_3)$ being naturally ordered.
}
\end{remark}

The following proposition will be an important ingredient in
the proof of Proposition~\ref{prop:hyper1}, while the last part
of its proof is used in the proof of
Theorem~\ref{thm:class}.

\begin{proposition}[The spinor equivariance]
\label{prop:planes}
There is a one-to-one correspondence between
oriented conformal structures on~$V_4$ and
definite $3$-planes $V_3$ in $V_6$ (with respect to~$Q$).
\end{proposition}

\begin{proof}
Given an oriented conformal structure on~$V_4$, define $V_3$ as
the corresponding space of self-dual $2$-forms on~$V_4$.

For the converse, we recall some quaternionic linear algebra. Under the
identification of the purely imaginary quaternions with~$\R^3$,
any element $\phi$ of $\SO (3)$ can be written as quaternionic conjugation
$\phi =\phi_u$,
\[ \R^3\ni x\longmapsto \phi_u(x)=ux\ou,\]
with some unit quaternion $u\in S^3\subset \HH$. The map
$u\mapsto \phi_u$ is the standard double covering $S^3\rightarrow\SO (3)$.

Given a definite $3$-plane in $V_6$, choose a naturally ordered
conformal symplectic triple $(A_1,A_2,A_3)$ spanning it. An orientation
on $V_4$ is then given by $A_i\wedge A_i$.

Let $q$ be a quaternionic coordinate for $V_4$ as in
Proposition~\ref{prop:pointmodel}, and $g=dx_1^2+dx_2^2+dx_3^2+dx_4^2$
the inner product on $V_4$ determined by $(A_1,A_2,A_3)$.

Suppose that $A_1',A_2',A_3'$ is another naturally ordered conformal
symplectic triple spanning the same $3$-plane $V_3$ and satisfying
$A_i'\wedge A_i'=A_i\wedge A_i$. Notice that $Q$ defines an inner
product on $V_3$ for which $(A_1,A_2,A_3)$ and $(A_1',A_2',A_3')$
are orthogonal bases consisting of vectors of equal length, and
defining the same orientation. Hence there is an element
$\phi_u\in\SO (3)$ such that
\[ \left( \begin{array}{c}A_1'\\A_2'\\A_3'\end{array} \right) =
\phi_u \left( \begin{array}{c}A_1\\A_2\\A_3\end{array} \right) .\]
By the preceding discussion this can be written as
\begin{eqnarray*}
\bi A_1'+\bj A_2'+\bk A_3' & = & u(\bi A_1+\bj A_2+\bk A_3)\ou\\
    & = & -\frac{1}{2}\, u\, dq\wedge d\oq\, \ou\\
    & = & -\frac{1}{2}\, d(uq)\wedge d(\overline{uq}).
\end{eqnarray*}
So a quaternionic coordinate corresponding to $(A_1',A_2',A_3')$
is given by~$uq$, which gives rise to the same inner product $g$
on $V_4$ since left multiplication on $\HH$ by a unit quaternion is
an isometry for $dx_1^2+dx_2^2+dx_3^2+dx_4^2$.

If the conformal symplectic triple $(A_1,A_2,A_3)$ is replaced by
$(vA_1,vA_2,vA_3)$, $v\in \R^+$, then the induced inner
product changes to~$vg$. This proves the proposition.
\end{proof}

\begin{proof}[Proof of Proposition~\ref{prop:hyper1}]
A (taut) contact sphere $(\alpha_1,\alpha_2,\alpha_3)$ on $U$
gives rise to a (conformal) symplectic triple
$(\Omega_1,\Omega_2,\Omega_3)$ on $U\times\R$ as described at the beginning of
this section. The complex structures $J_1,J_2,J_3$ defined on each
tangent space $T_x(U\times\R )$ depend smoothly
on the point~$x$, and thus define almost
complex structures on $U\times \R$, which are integrable
in the taut case, see~\cite{gego95}. So the statement concerning
contact spheres and conformal structures is immediate from the
preceding proposition. Notice that the orientation and conformal
class of the metric on $U\times \R$ are completely
characterised as the unique ones for which the $\Omega_i$ are
self-dual.

The statement about taut contact spheres and hyperk\"ahler structures
follows similarly; see~\cite{geig96} for an explicit argument.
A {\em non-taut} contact sphere determines at each point a
linear $2$-sphere worth
of complex structures, since pointwise the symplectic triple
$(\Omega_1,\Omega_2,\Omega_3)$ can be replaced by a conformal
symplectic triple spanning the same $3$-plane of skew-symmetric
forms; this does not change the space of corresponding almost
complex structures. Although this replacement can be done globally, leading
again to triples $(J_1,J_2,J_3)$ satisfying the quaternionic
identities, there does not seem to be a canonical choice for
doing it.

If two contact spheres are related by multiplication by the
function $v\co M\rightarrow \R^+$, the induced structures
on $U\times \R$ are related by the diffeomorphism
given by $(p,t)\mapsto (p, t+\log v(p))$.
\end{proof}

The following proposition explains how to go back from a
hyperk\"ahler structure and a tri-Liouville vector field
to a taut contact sphere.

\begin{proposition}
\label{prop:hyper2}
Let $(\Omega_1,\Omega_2,\Omega_3)$ be a
hyperk\"ahler structure on a $4$-manifold $W$ with
a nowhere zero tri-Liouville vector field~$Y$.
Then the equations
$\alpha_i=Y\ip\Omega_i$ define a naturally ordered taut
contact sphere $(\alpha_1,\alpha_2,\alpha_3)$ on any transversal
of~$Y$. Shifting the points of the transversal along the
orbits of $Y$ will change the taut contact sphere within
its conformal class.\qed
\end{proposition}

The proof is a straightforward computation, using $\Omega_i\wedge\Omega_i=
\Omega_j\wedge\Omega_j\neq 0$ and $\Omega_i\wedge\Omega_j=0$ for
$i\neq j$, cf.~\cite{geig96}, as well as the identities
$Y\ip\Omega_i^2=2\alpha_i\wedge d\alpha_i$, cf.\
identity~(\ref{eqn:cubic}). If the
hyperk\"ahler structure comes from a taut contact sphere,
then the above construction with $Y=\partial_t$ recovers that
contact sphere.

\section{Classification of taut contact spheres}
\label{section:class}
This section is largely devoted to the proof of
Theorem~\ref{thm:class}. As promised, this includes
a new proof that the universal cover of a closed
$3$-manifold carrying a taut contact sphere is diffeomorphic
to~$S^3$.

The strategy for the proof is as follows. From Theorem~\ref{thm:flat}
we know that a taut contact sphere on a closed $3$-manifold $M$
induces a flat metric on $M\times\R$.
With the help of the developing map of this
metric we are led to study hyperk\"ahler structures on Euclidean
$4$-space~$\E^4$, with the hyperk\"ahler metric being the flat Euclidean
metric. Parallel $2$-forms for this metric have constant coefficients
in Euclidean coordinates.
The classification of taut contact spheres is thus
reduced to a straightforward problem of determining
the possible tri-Liouville vector fields.

\begin{proof}[Proof of Theorem~\ref{thm:class}]
Let $M$ be a closed $3$-manifold
with a taut contact sphere, and let $(\alpha_1,\alpha_2,\alpha_3)$ be
the lifted taut contact sphere on its universal cover~$\widetilde{M}$.
Write $g$ for the induced metric on
$\widetilde{M}\times \R$. The fact that $g$ is flat and
$\widetilde{M}\times \R$
simply-connected implies that we have a developing map
$\Phi\co \widetilde{M}\times \R \rightarrow \E^4$
for this metric which is a local isometry.

If $W\subset \widetilde{M}\times \R$ is a domain on which $\Phi$ restricts
to a diffeomorphism, then $(\Phi |W)_*\partial_t$ is a vector field
$Y_W$ on the domain $\Phi (W)\subset \E^4$ generating a
$1$-parameter group of homotheties of the Euclidean metric (because
$L_{\partial_t}g=g$ by the construction of~$g$). Since homothetic
transformations of a Riemannian manifold are affine transformations
(this is easy to see for the Euclidean metric), $Y_W$ is the restriction
$Y|_{\Phi (W)}$ of a homothetic vector field~$Y$ defined on all
of~$\E^4$. Then $\partial_t$ and $\Phi^*Y$ are homothetic
vector fields for $g$ that coincide on the open set~$W$, which forces
$\partial_t= \Phi^*Y$ on all of $\widetilde{M}\times \R$.

A homothetic vector field on $\E^4$ vanishes at a single point,
and without loss of generality we may assume that $Y$ vanishes at~$0$.
Let $\pi\co \E^4\setminus\{ 0\} \rightarrow S^3_{\E}$ be the
projection onto the orbit space of~$Y$ (which we can identify with the
unit sphere $S^3_{\E}\subset \E^4$, for~$Y$ is a
genuinely expanding homothetic vector field and hence transverse to
any sphere centred at~$0$). Then the composition
\[ \widetilde{M}\times\{ 0\}\stackrel{\Phi}{\longrightarrow} \E^4\setminus
\{ 0\}\stackrel{\pi}{\longrightarrow} S^3_{\E} \]
is a local diffeomorphism and therefore,
$\widetilde{M}$ being simply-connected,
a diffeomorphism. So we have proved $\widetilde{M}\cong S^3$. Moreover,
the property $\partial_t=\Phi^*Y$ implies that $\Phi$
is a diffeomorphism from $S^3\times \R$ to $\E^4\setminus\{
0\}$, hence a global isometry.

From now on we write $S^3$ instead of~$\widetilde{M}$.
Now $\Phi$ sends $S^3\times\{ 0\}$ to some transversal of $Y$
in~$\E^4$, and by Proposition~\ref{prop:hyper2} the original contact
sphere is equivalent to the one induced on $S^3_{\E}$.

To simplify notation, we continue to write $e^t\alpha_i$, $\Omega_i$,
$J_i$ for the push-forwards of these objects to $\E^4\setminus\{ 0\}$,
and we identify $\partial_t$ with~$Y$. Thus $(J_1,J_2,J_3)$ defines
a hyperk\"ahler structure with respect to the Euclidean metric
$g_{\E}$, and $(\Omega_1,\Omega_2,\Omega_3)$ are the corresponding
K\"ahler forms. In particular, the $J_i$ and $\Omega_i$ are parallel
with respect to~$g_{\E}$, and thus have constant coefficients in
any linear coordinate system for~$\E^4$. As a consequence, there
are linear coordinates $x_0,x_1,x_2,x_3$ on $\E^4$ with respect to
which the formulae of Proposition~\ref{prop:pointmodel} hold
(with $A_i$ replaced by~$\Omega_i$). Observe that this forces
$x_0,x_1,x_2,x_3$ to be an orthonormal coordinate system with respect
to~$g_{\E}$.

Write $\psi_t$ for the flow of~$Y$. This flow commutes with the $J_i$ and
satisfies $\psi_t^*g_{\E}=e^tg_{\E}$ and $\psi_t^*
\Omega_i=e^t\Omega_i$, in particular $\psi_t^*(dz_1\wedge dz_2)=
e^t\, dz_1\wedge dz_2$.

Since the flow of $Y$ preserves $J_1$, we have a holomorphic vector field
$Y_{\C}$ on $\E^4=\C^2$ with $Y=2\,\mbox{\rm Re}
(Y_{\C})$. As a homothetic vector field vanishing at zero, $Y$
can be represented as a linear map, and this implies that $Y_{\C}$
corresponds to a complex linear map $\bz\mapsto\bY_{\C}\bz$ with
respect to the coordinate $\bz =\left( \begin{array}{c}z_1\\z_2
\end{array}\right)$.

A straightforward calculation shows that the condition $\psi_t^*g_{\E}
=e^tg_{\E}$ translates into the matrix $\bY_{\C}$ being of the form
\[ \bY_{\C}=\left( \begin{array}{cc}1/2 &\; 0\\0&\; 1/2\end{array}
\right) +\bZ_{\C}\]
with $\bZ_{\C}$ a skew-Hermitian matrix; the condition
$\psi_t^*(dz_1\wedge dz_2)=e^t\, dz_1\wedge dz_2$ forces $\bZ_{\C}$
to have zero trace. After a special unitary change of coordinates (which
does not change the expressions for the~$\Omega_i$) we may assume
that $\bZ_{\C}$ is in diagonal form, i.e.
\[ \bZ_{\C}=\left( \begin{array}{cc}\bi\nu &\; 0\\0&\; -\bi\nu\end{array}
\right) \]
with $\nu\in\R$. In the usual notation for vector fields
this means
\[ Y_{\C}=\Bigl( \frac{1}{2}+\bi\nu\Bigr) z_1\partial_{z_1}+
\Bigl( \frac{1}{2}-\bi\nu\Bigr) z_2\partial_{z_2}.\]
We conclude
\begin{eqnarray*}
e^t\alpha_1 & = & \partial_t\ip\Omega_1\\
     & = & Y\ip\bigl( \frac{\bi}{2} (dz_1\wedge d\oz_1+dz_2\wedge d\oz_2
               )\bigr)\\
     & = & \frac{\bi}{4}(z_1\, d\oz_1-\oz_1\, dz_1+z_2\, d\oz_2-\oz_2\, dz_2)
            -\frac{\nu}{2}\, d\bigl( |z_1|^2-|z_2|^2\bigr),\\
e^t(\alpha_2+\bi\alpha_3) & = & \partial_t\ip (\Omega_2+\bi\Omega_3)\\
       & = & Y\ip (dz_1\wedge dz_2)\\
       & = & \Bigl( \frac{1}{2}+\bi\nu \Bigr) z_1\, dz_2-
            \Bigl( \frac{1}{2}-\bi\nu \Bigr) z_2\, dz_1.
\end{eqnarray*}

\begin{remark}
{\rm
This analysis can be carried out locally.
As a result, the quaternionic formula in Theorem~\ref{thm:class} is
a universal local model for (naturally ordered) taut contact
spheres inducing a {\em flat\/} hyperk\"ahler metric.
}
\end{remark}

To translate the preceding equations into quaternionic notation, we observe
\[ q\cdot d\oq = z_1\, d\oz_1+z_2\, d\oz_2-z_1\, dz_2\bj +z_2\, dz_1\bj\]
and
\[ dq\cdot\oq = \oz_1\, dz_1+\oz_2\, dz_2+z_1\, dz_2\bj-z_2\, dz_1\bj\]
(using $\bj z=\oz\bj$ for $z\in\C$ and $\overline{z_1+z_2\bj}
=\oz_1-z_2\bj$), hence
\[ dq\cdot\oq -q\cdot d\oq = \oz_1\, dz_1-z_1\, d\oz_1+\oz_2\, dz_2-
z_2\, d\oz_2
     + 2(z_1\, dz_2-z_2\, dz_1)\bj .\]
We also observe
\[ q\bi\oq=\bigl( |z_1|^2-|z_2|^2\bigr) \bi-2\bi z_1z_2\bj .\]
Putting all this together, we find
\begin{eqnarray*}
e^t(\bi\alpha_1+\bj\alpha_2+\bk\alpha_3) & = &
                   e^t(\bi\alpha_1+(\alpha_2+\bi\alpha_3)\bj)\\
            & = & \frac{1}{4}(dq\cdot\oq-q\cdot d\oq)
                   -\frac{\nu}{2} \,d(q\bi\oq ).
\end{eqnarray*}
This shows that, up to conformal equivalence and diffeomorphism, the taut
contact spheres on $S^3$ are as described in Theorem~\ref{thm:class}.
The fact that different non-negative values of $\nu$ give non-isomorphic
contact spheres follows from the corresponding classification of
taut contact circles in~\cite{gego95}, \cite{gego02}. We shall
be a bit more explicit about this point below, where we give
a pictorial description of the moduli spaces in
question. This will include a synthetic method for
determining the modulus, independent of our previous papers. 

Next we want to show that $\nu$ and $-\nu$ correspond
to equivalent structures.
The diffeomorphism of $S^3$ given by
$q\mapsto q\bj$ is isotopic to the identity and pulls the
quaternionic $1$-form with parameter value $\nu$ to that with
value $-\nu$, because $\bj$ anticommutes with $\bi$.

After having determined the possible lifted taut contact
spheres on $\widetilde{M}\cong S^3$, we now consider the taut contact
sphere on $M$ itself. From Theorem~\ref{thm:classold} we know
that $M$ has to be a left-quotient of~$\SU (2)$.
The induced flat hyperk\"ahler structure on $M\times \R$ lifts to
just such a structure on $S^3\times \R$, invariant under the
deck transformation group~$\Gamma$. As was already argued
in~\cite{gego95} for taut contact circles, this implies that in the complex
coordinates $(z_1,z_2)$ which give a normal form as described above,
one has $\Gamma\subset\SU (2)$. Moreover, again as in~\cite{gego95},
the parameter $\nu$ is forced to be zero for non-abelian~$\Gamma$,
and it can take any value for the group $\Gamma =\Z_m\subset\SU (2)$
generated by
\[ \left( \begin{array}{cc}\varepsilon &\; 0\\ 0&\;\varepsilon^{-1} \end{array}
\right) ,\]
where $\varepsilon$ is some $m$th root of unity.

\vspace{2mm}

Now we address the homogeneity issue in the statement of
Theorem~\ref{thm:class}.
For any unit quaternion $u$ we have
\[ u\Bigl( \frac{1}{2}(dq\cdot \oq -q\cdot d\oq )-\nu\, d (q\bi\oq )
\Bigr) \ou =\frac{1}{2}(d(uq)\cdot\overline{uq}-uq\cdot d(
\overline{uq}))-\nu\, d(uq\bi\overline{uq}).\]
Therefore, if two taut contact spheres
are related by an element $\phi_u$ of $\SO (3)$
as in the proof of Proposition~\ref{prop:planes}, then one
is the pull-back of the other under a diffeomorphism of
$\Gamma\backslash\SU (2)$ induced by the map
$q\mapsto uq$. In particular, this
$\SO (3)$-action on taut contact spheres shows that any taut contact
sphere can be swept out by great circles, all defining isomorphic taut
contact circles.

This concludes the proof of Theorem~\ref{thm:class}.
\end{proof}

This theorem allows us to define a map from the space of taut contact
spheres to that of taut contact circles: simply pick any great circle.

\begin{theorem}
\label{thm:moduli}
Let $M$ be any left-quotient of $\SU (2)$.
The map sending a taut contact sphere on $M$ to any of its great
circles induces an embedding of the moduli space of taut contact spheres
on $M$ into the moduli space of taut contact circles.
\end{theorem}

\begin{proof}
Recall the classification (up to diffeomorphism and homothety, i.e.\
conformal equivalence and rotation)
of taut contact circles $(\alpha_2,\alpha_3)$
on left-quotients of~$S^3$,
see \cite{gego95}, \cite{gego02}. On the lens spaces $L(m,m-1)$ we have the
continuous family of taut contact circles induced by the
$\Z_m$-invariant complex $1$-form
\[ \alpha_2+\bi\alpha_3=\Bigl( \frac{1}{2} +\delta \Bigr)
z_1\, dz_2- \Bigl( \frac{1}{2} -\delta \Bigr) z_2\, dz_1\]
(restricted to $S^3\subset\C^2$) with $\delta\in \C$,
$-1/2<\mbox{\rm Re}(\delta )<1/2$, modulo replacing $\delta$ by $-\delta$,
which corresponds to the diffeomorphism defined by $(z_1,z_2)\mapsto
(z_2,z_1)$ and changing from $(\alpha_2,\alpha_3)$ to $(-\alpha_2,-\alpha_3)$.
So from Theorem~\ref{thm:class} we see that
\begin{itemize}
\item only the taut contact
circles with $\delta$ purely imaginary extend to taut contact spheres,
\item this extension is unique up to automorphisms of $\alpha_2+
\bi\alpha_3$, and
\item two taut contact circles giving rise to isomorphic taut contact
spheres must be isomorphic.
\end{itemize}
This proves the theorem in the given case. For the
abelian left-quotients the map on moduli
spaces amounts to the inclusion
\[
\begin{array}{rcl}
\R_0^+ & \longrightarrow & \{\delta\in\C\,\co
     -1/2<\mbox{\rm Re}(\delta )<1/2\}/_{\delta\sim -\delta}\\
\nu & \longmapsto & \bi\nu.
\end{array} \]

The moduli space of taut contact circles
on $L(m,m-1)$ also includes a discrete family described by
\[ \alpha_2+\bi\alpha_3=nz_1\, dz_2-z_2\, dz_1+z_2^n\, dz_2,\]
where $n$ ranges over the natural numbers congruent $-1$ mod~$m$.
These contact circles, however, do not extend to any taut contact sphere,
so they are not being considered here.
\end{proof}

\begin{figure}[h]
\centering
\includegraphics[scale=0.35]{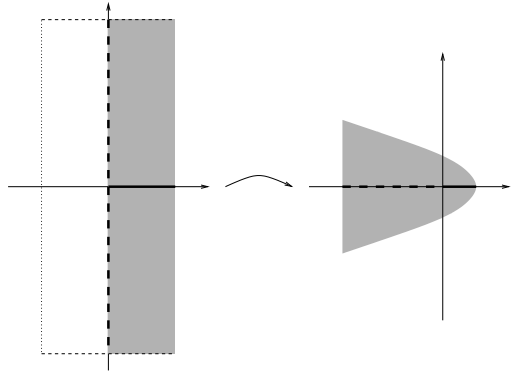}
  \caption{The moduli spaces of taut contact circles and spheres.}
  \label{figure:moduli}
\end{figure}

The moduli spaces described in the foregoing proof
(without the discrete family) are
illustrated on the left-hand side of Figure~\ref{figure:moduli}.
The moduli space of taut contact circles is the orbifold
\[ \{\delta\in\C\,\co -1/2<\mbox{\rm Re}(\delta )<1/2\}/_{\delta\sim -\delta}.\]
The half-line
\[ \{\delta\in\bi\R\}/_{\delta\sim -\delta}\]
constituting the moduli space of taut contact spheres is shown
as a dashed line. The real part
\[ \{\delta\in\R\,\co -1/2<\delta <1/2\}/_{\delta\sim -\delta}\]
of this moduli space, shown in bold,
corresponds to so-called Cartan structures. The origin
represents the unique $3$-Sasakian structure. These structures will
be discussed in Section~\ref{section:Sasakian}.

Under the mapping $\delta\mapsto\delta^2$, the moduli space of taut
contact circles is mapped bijectively to the interior of
the parabola
\[ \Bigl\{ x+\bi y\in\C\,\co x=\frac{1}{4}-y^2\Bigr\}\]
shown on the right-hand side of Figure~\ref{figure:moduli}. So this
singular mapping flattens the cone point and yields
a representation of the moduli space as an open subset of~$\C$.

Here is the promised synthetic characterisation of the modulus $\nu$ of
a taut contact sphere on~$S^3$. One can define a
{\em canonical slice} inside the hyperk\"ahler manifold $S^3\times\R$
as the subset $\{ \|\partial_t\|=1\}$. A straightforward calculation shows
that for the taut contact sphere of modulus $\pm\nu$ we have
\[ \|\partial_t\|^2=\Bigl(\frac{1}{4}+\nu^2\Bigr)\cdot\bigl(
|z_1|^2+|z_2|^2\bigr) .\]
So the canonical slice is isometric with
the $3$-sphere of radius $2/\sqrt{1+4\nu^2}$.
The canonical slice of maximal radius corresponds to the unique
$3$-Sasakian structure.

The canonical slice has the property that the $1$-forms
$e^t\alpha_i$ induce on it a $1$-normalised taut contact sphere.
We conclude that the family in Theorem~\ref{thm:class} induces on
any Euclidean sphere a $c$-normalised taut contact sphere, for
some constant~$c$.

Recall now Definition \ref{defn:short-long}. Once the contact
sphere is $1$-normalised, the canonical slice is given by $\{ t=0\}$
and we see that the hyperk\"ahler metric induces the long
metric $g_l$ on it. Thus, on a
closed $3$-manifold the long metric is always spherical.

The short metric can be written as $g_s=g_l-\beta^2$. We claim that
$\beta$ is invariant under the maps $q\mapsto uq$ with $u\in
S^3$, which shows that $g_s$ is a {\em Berger metric.} To see
this invariance property, notice that the system of structure equations
(\ref{eqn:structure}) is invariant under rotations of the triple
$(\alpha_1,\alpha_2,\alpha_3)$, and we have seen that left
multiplication by unit quaternions $u$ induces such rotations.

The contact sphere with $\nu =0$ in Theorem~\ref{thm:class} is
invariant under right multiplication by unit quaternions, which
makes the corresponding $\beta$ bi-invariant, but recall that in
this case $\beta\equiv 0$ and so there is no contradiction here.

\vspace{2mm}

We briefly expand on the point in the proof of Theorem~\ref{thm:moduli}
concerning the
automorphisms of $\alpha_2+\bi\alpha_3$. From the arguments in
\cite[Section~5.2]{gego95} it follows that automorphisms $\phi$
of $\alpha_2+\bi\alpha_3$ are $\C$-linear maps of
$\C^2$, in fact, elements of $\mathrm{SL}_2\C$.
Given an $\alpha_1$ extending
$(\alpha_2,\alpha_3)$ to a taut contact sphere, the other extensions
are given by $\phi^*\alpha_1$.

For $S^3$ and $\nu =0$, any element of $\mathrm{SL}_2\C$
defines an automorphism of $\alpha_2+\bi\alpha_3$, so the
possible extensions are parametrised by $\mathrm{SL}_2\C/\SU (2)$.

For $S^3$ and $\nu\neq 0$, the condition that $\phi$ preserve
$\alpha_2+\bi\alpha_3$ forces it to be a diagonal map
$\phi_c(z_1,z_2)=(cz_1,c^{-1}z_2)$.
The same is true for the lens spaces $\Z_m\backslash\SU (2)$
other than $S^3$, even in the case $\nu =0$; here the condition for $\phi$
to be diagonal follows from the fact that it has to lie in the normaliser
of~$\Z_m$.
In these two cases, the possible extensions are
parametrised by $\R^+$, since
\[ \phi_c^*\alpha_1=
\frac{\bi}{4}\bigl( |c|^2(z_1\, d\oz_1-\oz_1\, dz_1)+|c|^{-2}
(z_2\, d\oz_2-\oz_2\, dz_2)
             \bigr)
            -\frac{\nu}{2}\, d\bigl( |c|^2|z_1|^2-|c|^{-2}|z_2|^2\bigr) .
\]

Finally, for the quotients of $\Gamma\backslash\SU (2)$ with
$\Gamma$ non-abelian, the fact that $\phi$ lies in the normaliser
of $\Gamma$ forces it to be plus or minus the identity map. Thus the
extension $\alpha_1$ is unique.

\section{Non-flat metrics}
\label{section:nonflat}
In this section we describe two constructions of taut contact spheres
giving rise to non-flat hyperk\"ahler metrics. These examples constitute
a proof of Theorem~\ref{thm:nonflat}, since a non-flat
hyperk\"ahler metric is not even conformally flat
(cf.\ the proof of Theorem~\ref{thm:flat}).

\subsection{The Helmholtz equation on the $2$-sphere}
\label{section:Helmholtz}
We have seen in Proposition~\ref{prop:hyper2} how to
construct a taut contact sphere corresponding to a suitable
hyperk\"ahler structure $(\Omega_1,\Omega_2,\Omega_3)$. Notice
that in terms of the holomorphic structure given by~$J_1$, an
equivalent description of this hyperk\"ahler structure is given by
a holomorphic symplectic form $\Omega=\Omega_2+\bi \Omega_3$ and a
closed real $(1,1)$-form $\Omega_1$ with $2\Omega_1^2=\Omega\wedge
\overline{\Omega}$.

We now make the following ansatz: Let $z_1,z_2$ be complex
coordinates on $\C^2$ and identify $\C^2$ with
$\R^3 \times \R$ by equating the $\R$-direction
with the real part of~$2z_1$, so that
$\partial_t=\mbox{\rm Re}(\partial_{z_1})$ and
$\psi_t(z_1,z_2)=(z_1+t/2,z_2)$. Set
\[ \Omega=\Omega_2+\bi\Omega_3=\lambda (z_1,z_2)\, dz_1\wedge dz_2\]
with $\lambda (z_1,z_2)$ a nowhere zero holomorphic function, and
\[ \Omega_1=\frac{\bi}{2} \partial\opartial H(z_1,z_2)\]
with $H(z_1,z_2)$ a real-valued function.

To satisfy the conditions $\psi_t^*\Omega_i=e^t\Omega_i$ it is sufficient
to have
\begin{equation}
\label{eqn:homogeneity}
\psi_t^*\lambda =e^t\lambda,\;\;\; \psi_t^*H=e^tH.
\end{equation}
Further, the identity $2\Omega_1^2=\Omega\wedge\overline{\Omega}$ is
equivalent to the complex Monge-Amp\`ere equation
\begin{equation}
\label{eqn:MA}
\det \left( \begin{array}{cc}
H_{z_1\oz_1} & \; H_{z_1\oz_2}\\
H_{z_2\oz_1} & \; H_{z_2\oz_2}
\end{array} \right) =\lambda\overline{\lambda}.
\end{equation}
If all these conditions are met, Proposition~\ref{prop:hyper2}
tells us how to recover the corresponding taut contact sphere.

We satisfy condition (\ref{eqn:homogeneity}) by taking $\lambda =e^{2z_1}$
and simplifying the ansatz further to
\begin{eqnarray*}
H(z_1,z_2) & = & e^{2\mbox{\scriptsize\rm Re}(z_1)}h(2\, \mbox{\rm Im}(z_1),
                    2\,\mbox{\rm Re}(z_2))\\
   & = & e^{z_1+\oz_1}h(\bi\oz_1-\bi z_1,z_2+\oz_2),
\end{eqnarray*}
with $h(s_1,s_2)$ a function of two real variables. Then, writing
$h_i$ for the partial derivatives $h_{s_i}=\partial h/\partial_{s_i}$ etc.,
we have
\begin{eqnarray*}
\partial\opartial H & = & e^{z_1+\oz_1}\bigl( (h+h_{11})\, dz_1\wedge d\oz_1+
                          (h_2-\bi h_{12})\, dz_1\wedge d\oz_2\\
   & & \;\;\;\;\;\;\;\;\;\mbox{} +(h_2+\bi h_{12})\,
                      dz_2\wedge d\oz_1+h_{22}\, dz_2
                     \wedge d\oz_2\bigr) ,
\end{eqnarray*}
and equation (\ref{eqn:MA}) becomes
\begin{equation}
\label{eqn:MA2}
1=\left| \begin{array}{cc}
h+ h_{11}      & \; h_2-\bi h_{12}\\
h_2+\bi h_{12} & \; h_{22}
\end{array} \right|
= \left| \begin{array}{cc}
h_{11} & \; h_{12}\\
h_{12} & \; h_{22}
\end{array} \right|
+\left| \begin{array}{cc}
h   & \; h_2\\
h_2 & \; h_{22}
\end{array} \right| .
\end{equation}

The reader may check directly that the function
\[ h(s_1,s_2)=\cos s_1\, \int_{1}^{s_2/\cos s_1} \sqrt{\xi^2-1}\, d\xi \]
satisfies this equation. In the appendix we shall explain
how to derive this sample solution, and verify that it leads to a non-flat
hyperk\"ahler metric. We arrive at this example
after proving a general result that relates solutions of
(\ref{eqn:MA2}) to solutions of the Helmholtz equation on the $2$-sphere.
\subsection{The Gibbons-Hawking ansatz}
\label{section:GH}
The Gibbons-Hawking ansatz \cite{giha78} starts with a
positive function $V$ in three variables $x_1,x_2,x_3$
(locally on~$\R^3$)
and a triple of functions
$b_1,b_2,b_3$ in the same variables, satisfying the
condition
\begin{equation}
\label{eqn:curl}
\nabla V=-\mbox{\rm curl} (b_1,b_2,b_3),
\end{equation}
so that in particular
\[ \Delta V=\mbox{\rm div}(\nabla V)=-\mbox{\rm div}(
\mbox{\rm curl} (b_1,b_2,b_3))=0,\]
i.e.\ $V$ is harmonic.
Set $\beta =b_1\, dx_1+b_2\, dx_2+b_3\, dx_3$ and consider
the triple of $2$-forms on $\R^3\times\R$ (with $\theta$
denoting the $\R$-coordinate) defined by
\[ \Omega_i=(d\theta +\beta)\wedge dx_i+V\, dx_j\wedge dx_k.\]
The relation between $V$ and $\beta$ implies that these $2$-forms
are closed. Writing
\[ \theta_0:= V^{-1/2}(d\theta +\beta )\;\;\mbox{\rm and}\;\;
\theta_i:=V^{1/2}dx_i,\;\; i=1,2,3,\]
we have
\[ \Omega_i=\theta_0\wedge\theta_i+\theta_j\wedge\theta_k,\]
which shows that $(\Omega_1,\Omega_2,\Omega_3)$ is a conformal
symplectic triple. The corresponding hyperk\"ahler metric is
\[ \theta_0^2+\theta_1^2+\theta_2^2+\theta_3^2=
V^{-1}(d\theta +\beta )^2+V(dx_1^2+dx_2^2+dx_3^2).\]
Observe the formal similarity with the metric in~(\ref{eqn:metric}).

The translational invariance of this metric in $\theta$-direction is
obvious.\footnote{In particular, the metric descends to
$\R^3\times S^1$; this is the usual form of
the Gibbons-Hawking ansatz.} A homothety can be built into this ansatz by
choosing $\beta$ appropriately.

Here is an example. As domain $U_0\subset\R^3$ we take
the half-space given by the condition $x_1>0$, and $V(x_1,x_2,x_3):=
x_1$ is our positive function on that domain. Set $(b_1,b_2,
b_3)=(0,x_3,0)$, that is, $\beta =x_3\, dx_2$. Then
\[ \nabla V= (1,0,0)=-\mbox{\rm curl} (b_1,b_2,b_3),\]
i.e.\ condition~(\ref{eqn:curl}) is satisfied. Then the Gibbons-Hawking ansatz
yields the hyperk\"ahler structure\footnote{The change in notation from
$\theta$ to $x_0$ is meant to emphasise that we want this to be
an $\R$-, not an $S^1$-coordinate.}
\begin{eqnarray*}
\Omega_1 & = & (dx_0+x_3\, dx_2)\wedge dx_1+x_1\, dx_2\wedge dx_3,\\
\Omega_2 & = & dx_0\wedge dx_2+x_1\, dx_3\wedge dx_1,\\
\Omega_3 & = & (dx_0+x_3\, dx_2)\wedge dx_3+x_1\, dx_1\wedge dx_2,
\end{eqnarray*}
with corresponding hyperk\"ahler metric
\[ g=\frac{1}{x_1}(dx_0 +x_3\, dx_2)^2+x_1\, (dx_1^2+dx_2^2+dx_3^2).\]
A tri-Liouville vector field for this
hyperk\"ahler structure is
\[ Y:=\frac{2}{3}x_0\partial_{x_0}+\frac{1}{3}(x_1\partial_{x_1}+
x_2\partial_{x_2}+x_3\partial_{x_3}).\]
By Proposition~\ref{prop:hyper2}, the equations $\alpha_i=Y\ip\Omega_i$,
$i=1,2,3$, define a taut contact sphere on any transversal $U$ to~$Y$.
By construction, $g$ is in turn the hyperk\"ahler metric
on $U\times\R$ induced by this taut contact sphere, with
the $\R$-factor now corresponding to the flow lines of~$Y$.

The surface $\Sigma:=\{ x_2=x_3=0\}\subset U_0\times\R$ is
totally geodesic for the metric~$g$, since it is the fixed point set of the
isometric involution $(x_0,x_1,x_2,x_3)\mapsto (x_0,x_1,-x_2,-x_3)$.
The metric on $\Sigma$ induced by $g$ is
\[ \frac{1}{x_1}\, dx_0^2+x_1\, dx_1^2,\]
and, from the well-known formula for computing the Gau{\ss} curvature
of a metric in diagonal form, one obtains $K_{\Sigma}=-1/x_1^3$.
This proves that $g$ is non-flat.

\subsection{Comparison of the two examples}
We observe that our ansatz from
Section~\ref{section:Helmholtz} can also be related to
a harmonic function on~$\R^3$:

Write $\Delta^{S^2}$ for the Laplacian on the unit $2$-sphere
in~$\R^3$, and $\Delta^{\R^3}$ for the Laplacian
on~$\R^3$. Moreover, let $\rho$ be the radial coordinate
on~$\R^3$. Given a differentiable function $\tu\co
\R^3\rightarrow \R$, we have the relation
\[ \Delta^{S^2}(\tu |_{S^2})=(\Delta^{\R^3}\tu-\tu_{\rho\rho}
-2\tu_{\rho})|_{S^2}; \]
the last summand reflects the fact that both principal curvatures
of~$S^2$ are equal to~$1$.

If $\tu$ is of the form
\[ \tu (x_1,x_2,x_3)=\rho u(x_1/\rho ,x_2/\rho ,x_3/\rho ),\]
then this last equation simplifies to
\[ \Delta^{S^2}u+2u=(\Delta^{\R^3}\tu)|_{S^2}.\]
In that particular case, $\tu$ is homogeneous of degree~$1$
in $\rho$, hence $\Delta^{\R^3}\tu$ is homogeneous
of degree~$-1$. This implies that the vanishing of
$(\Delta^{\R^3}\tu)|_{S^2}$ is sufficient for the
vanishing of $\Delta^{\R^3}\tu$ on all~$\R^3$.
In conclusion, solutions to our Monge-Amp\`ere equation~$(\ref{eqn:MA2})$
are --- by the preceding discussion and Proposition~\ref{prop:helmholtz}
in the appendix ---
in direct correspondence with harmonic functions on $\R^3$
that are homogeneous of degree~$1$ in~$\rho$.

Finally, notice that both examples admit a tri-Hamiltonian
vector field: in the example of Section~\ref{section:Helmholtz},
this is the vector field $\partial_{x_3}$, with $x_3:=\mathrm{Im}(z_2)$;
in the example of Section~\ref{section:GH}, it is~$\partial_{x_2}$.
As we plan to show in a forthcoming paper, all taut contact
spheres with such a symmetry can be related to a Gibbons-Hawking
ansatz, although this relation, in the Helmholtz case, is far from
straightforward.
\section{A Bernstein problem}
\label{section:Bernstein}
Are there any $1$-normalised taut contact spheres --- on a suitable open
domain~$U$ --- inducing a non-flat metric $g$ on $U\times\R$ such
that the long metric $g_l$ (see Definition~\ref{defn:short-long})
induced on $U\equiv U\times\{ 0\}$ is complete?
This kind of completeness question is known as a {\it Bernstein problem},
see~\cite{cala70}.

In order to appreciate the difficulty of this question, it is helpful
to consider how the ansatz in Section~\ref{section:Helmholtz} has
to be extended if it is to be of any use in providing an answer.
First of all, we ensured the homogeneity of $\partial\opartial H$ by
taking $H$ to be homogeneous. However, the example
\[ H=z_1e^{\oz_1}+\oz_1e^{z_1},\]
where
\[ \opartial H=\bigl( z_1e^{\oz_1}+e^{z_1}\bigr)\, d\oz_1\;\;
\mbox{\rm and}\;\; \partial\opartial H=\bigl( e^{\oz_1}+e^{z_1}\bigr)\, dz_1
\wedge d\oz_1,\]
shows that it is perfectly possible for $\partial\opartial H$ to be
homogeneous in $e^{\mathrm{Re}(z_1)}$ without either $H$ or $\opartial H$
having this property. Secondly, we took $H$ to be independent of
$\mathrm{Im}(z_2)$, but to discuss the general case one needs to allow
the auxiliary function $h$ to depend on all three variables
$\mathrm{Im}(z_1)$, $\mathrm{Re}(z_2)$, and~$\mathrm{Im}(z_2)$. Thirdly,
while a potential $H$ always exists locally, it need not exist globally.

The question we just raised has a venerable history. In
\cite{cala75} Calabi describes a construction of
{\it K\"ahler-Einstein} metrics, and in particular Ricci-flat K\"ahler
metrics, on complex tubular domains $D\times \bi\R^n\subset
\C^n$, where $D$ is some connected, open subset of
$\R^n$. The construction is based on a real
Monge-Amp\`ere equation (with a constant on the right-hand side),
or equivalently, a complex Monge-Amp\`ere equation for a function
depending only on the real parts of $n$ complex variables. The
resulting metrics are invariant under the group of translations
along $\bi\R^n$. Earlier results of Ca\-labi~\cite{cala58}
allowed him to show that, in the Ricci-flat case, metrics obtained
via this construction can never be complete, except for the
trivial case with $D=\R^n$ and a flat metric.

Here is a partial answer to this Bernstein problem,
dealing with contact spheres that possess additional symmetries. In place of
the long metric $g_l$ we consider the short metric $g_s$.
Completeness of $g_s$ is a stronger condition than
completeness of~$g_l$. The following is part (a)
of Theorem~\ref{thm:Bernstein}.

\begin{theorem}
\label{thm:Bernstein-Killing}
Let $(\omega_1,\omega_2,\omega_3)$ be a $1$-normalised taut contact sphere
on a $3$-mani\-fold $U$ with the property that the short metric
$g_s :=\omega_1^2+
\omega_2^2+\omega_3^2$ is complete and admits a Killing
field $X\not\equiv 0$ that preserves each form, i.e.
\[ L_X\omega_1=L_X\omega_2=L_X\omega_3=0.\] 
Then $U$ is compact and hence a left-quotient of $\SU (2)$,
and $(\omega_1,\omega_2,\omega_3)$ is isomorphic to
one of the taut contact spheres described in Theorem~\ref{thm:class}.
\end{theorem}

The proof of this theorem will take up the remainder of this section.
We begin by observing that $X$ does not have any zeros, which
can be seen as follows. Arguing by contradiction, assume that
$p\in U$ is a point with $X(p)=0$. Then the flow of $X$ preserves
the distance spheres from~$p$, and is then necessarily a rotation
about an axis through~$p$ (in geodesic normal coordinates). This is
incompatible with the fact that this flow preserves the coframe
$(\alpha_1,\alpha_2,\alpha_3)$.

The fact that $X$ does not have any zeros means that
\[ \Lambda:=\bigl( \omega_1(X)^2+\omega_2(X)^2+\omega_3(X)^2\bigr)^{1/2}\]
defines a function $\Lambda\co U\rightarrow\R^+$.
Set
\[ \alpha_i:=\omega_i/\Lambda,\;\; i=1,2,3,\]
so that
\[ \alpha_1(X)^2+\alpha_2(X)^2+\alpha_3(X)^2\equiv 1.\]
Notice that the $\alpha_i$ satisfy the
equations~(\ref{eqn:structure}) with that very~$\Lambda$, and they are likewise
invariant under the flow of~$X$.

Consider the map
\[ \Psi:=(x_1,x_2,x_3):=\bigl(\alpha_1(X),\alpha_2(X),\alpha_3(X)\bigr)\co
U\longrightarrow S^2.\]
It is clear that the differential of $\Psi$ satisfies $T\Psi (X)=0$,
i.e.\ each flow line of $X$ is mapped to a single point in~$S^2$.
Write $\|.\|_s$ for the length of tangent vectors to $U$
with respect to the short metric~$g_s$, and $\|.\|_{S^2}$ for the length
of tangent vectors to $S^2\subset\R^3$
with respect to the standard metric $dx_1^2+dx_2^2+dx_3^2$.

\begin{lemma}
\label{lem:1}
If $Z\in T_pU$ is a tangent vector $g_s$-orthogonal to $X$, then
\[ \| T\Psi(Z)\|_{S^2}\geq\| Z\|_s.\]
\end{lemma}

\begin{proof}
The $X$-invariance of the $\alpha_i$ gives, with the Cartan formula for the
Lie derivative,
\[ dx_i=d(\alpha_i(X))=L_X\alpha_i-X\ip d\alpha_i=-X\ip d\alpha_i.\]
All the following computations are made at the single point~$p$.
Rotate the contact sphere so that at that point $p$ we have
$\alpha_1(X)=1$ and $\alpha_2(X)=\alpha_3(X)=0$, i.e.\
$\Psi (p)=(1,0,0)$. Then, with~(\ref{eqn:structure}),
\begin{eqnarray*}
dx_1 & = & -\beta (X)\alpha_1+\beta,\\
dx_2 & = & -\beta (X)\alpha_2+\Lambda\,\alpha_3,\\
dx_3 & = & -\beta (X)\alpha_3-\Lambda\,\alpha_2.
\end{eqnarray*}
The condition that $Z$ be $g_s$-orthogonal to $X$ means that
$\alpha_1(Z)=0$. Hence
\begin{eqnarray*}
dx_1(Z) & = & \beta (Z),\\
dx_2(Z) & = & -\beta (X)\alpha_2(Z)+\Lambda\,\alpha_3(Z),\\
dx_3(Z) & = & -\beta (X)\alpha_3(Z)-\Lambda\,\alpha_2(Z).
\end{eqnarray*}
From $\Psi (p)=(1,0,0)$ we have $dx_1(Z)=0$. Then
\begin{eqnarray*}
\| T\Psi (Z)\|_{S^2}^2 & = & dx_1(Z)^2+dx_2(Z)^2+dx_3(Z)^2\\
  & = & dx_2(Z)^2+dx_3(Z)^2\\
  & = & \bigl(\beta (X)^2+\Lambda^2\bigr)\cdot\bigl(\alpha_2(Z)^2
         +\alpha_3(Z)^2\bigr)\\
  & = & \bigl(\beta (X)^2+\Lambda^2\bigr)\cdot\bigl(\alpha_1(Z)^2
         +\alpha_2(Z)^2+\alpha_3(Z)^2\bigr)\\
  & \geq & \Lambda^2\bigl(\alpha_1(Z)^2
         +\alpha_2(Z)^2+\alpha_3(Z)^2\bigr)\\
  & = & \| Z\|_s^2.\;\;\qed
\end{eqnarray*}
\renewcommand{\qed}{}
\end{proof}

This lemma implies in particular that $T\Psi$ has full rank at
every point, so $\Psi (U)$ is an open subset of~$S^2$.

\begin{lemma}
The map $\Psi\co U\rightarrow S^2$ is surjective.
\end{lemma}

\begin{proof}
Since $\Psi (U)\subset S^2$ is open, it suffices to show
that $\Psi (U)$ is complete, i.e.\ that every path
$\gamma \co [0,1)\rightarrow \Psi (U)$ of finite length
has a limit point inside $\Psi (U)$.
Let $\tilde{\gamma}\co [0,t_0)\rightarrow U$ be a maximal lift
of such a path, $g_s$-orthogonal to $X$. By the previous lemma,
this lift is non-empty and of finite length in the short metric. Since
$U$ is complete, we deduce that $t_0=1$ and that the lift $\tilde{\gamma}$
has a limit point in $U$. Therefore $\gamma$ has a limit point in $\Psi (U)$.
\end{proof}

\begin{lemma}
For each $q\in S^2$, the preimage $\Psi^{-1}(q)\subset U$ is
a single orbit of~$X$.
\end{lemma}

\begin{proof}
Arguing by contradiction, we assume that $q\in S^2$ is a point
in the image of two distinct orbits ${\mathcal O}_0,{\mathcal O}_1$
of~$X$. Let $\tilde{\gamma}
\co [0,1]\rightarrow U$ be
a path joining these two orbits; $\gamma:=\Psi\circ\tilde{\gamma}
\co [0,1]\rightarrow S^2$ is then a loop based at~$q$. By the argument
in the proof of the preceding lemma we may assume that $\tilde{\gamma}$
is orthogonal to~$X$. Let $\gamma_s$, $s\in [0,1]$,
be a homotopy of $\gamma=\gamma_0$ rel~$\{ 0,1\}$ to the constant
path~$\gamma_1$ at~$q$, and $\tilde{\gamma}_s$ the corresponding homotopy
of lifts orthogonal to $X$ with initial point
$\tilde{\gamma}_s(0)=\tilde{\gamma}(0)$ for all $s\in [0,1]$.
The endpoints $\tilde{\gamma}_s(1)$ form a smooth path in~$\Psi^{-1}(q)$.

The completeness of $(U,g_s)$ entails that the flow of $X$ is complete,
and this in turn ensures that $\tilde{\gamma}_s(1)\in{\mathcal O}_1$
for all~$s$. But the constant path $\gamma_1$ lifts to the
constant path $\tilde{\gamma}_1$ at $\tilde{\gamma}(0)\in {\mathcal O}_0$,
hence $\tilde{\gamma}_1(1)\in{\mathcal O}_0\cap
{\mathcal O}_1$, which is impossible.
\end{proof}

As preimages of single points, the orbits of $X$ are closed subsets of~$U$.
This means that every orbit is either periodic or a proper embedding
of $\R$ in~$U$. If there is a periodic orbit ${\mathcal O}_0$,
all other orbits remain at a bounded distance from~${\mathcal O}_0$,
since the flow of $X$ is by isometries. This precludes proper embeddings
of $\R$, i.e.\ in this case all orbits must be periodic,
and $U$ is compact, as asserted in Theorem~\ref{thm:Bernstein-Killing}.

It remains to consider the complementary case, when all orbits are proper
embeddings of~$\R$.
In this case, the time-$1$ map of the flow of $X$ will disjoin any
sufficiently small compact set $K$ from itself, since each single orbit
through $K$ is proper and the flow is by isometries.
It follows that this time-$1$ map defines a free and
properly discontinuous $\Z$-action on $U$.

The quotient manifold under this action is, by the first case,
a compact manifold with a taut contact sphere, and hence
a left-quotient of $\SU (2)$. Such manifolds do not admit infinite
covers. In other words, this second case cannot occur.

This completes the proof of Theorem~\ref{thm:Bernstein-Killing}.

\section{$3$-Sasakian structures}
\label{section:Sasakian}
In this section we consider taut contact spheres
$(\alpha_1,\alpha_2,\alpha_3)$ on a $3$-dimensional domain $U$
that satisfy the stronger condition $\alpha_i\wedge d\alpha_j=0$
for $i\neq j$. The conditions for a naturally ordered taut contact
sphere then imply that $d\alpha_i=\Lambda\,\alpha_j\wedge\alpha_k$
for some $\Lambda\co U\rightarrow \R^+$. Differentiation
of this equation yields
\[ 0=d(d\alpha_i)=d\Lambda\wedge\alpha_j\wedge\alpha_k+\Lambda\,
d\alpha_j\wedge\alpha_k-\Lambda\,\alpha_j\wedge d\alpha_k
=d\Lambda\wedge\alpha_j\wedge\alpha_k.\]
It follows that $d\Lambda\equiv 0$, so $\Lambda\equiv c$ for some constant
$c$ and the contact sphere is $c$-normalised.

Recall that the {\bf Reeb vector field} $R$ of a contact
form $\alpha$ is defined by $\alpha (R)\equiv 1$ and $R\ip
d\alpha\equiv 0$. We now have the following simple lemma,
where $R_i$ denotes the Reeb vector field of~$\alpha_i$,
and $\beta$ is the $1$-form from equation~(\ref{eqn:structure}).
The proof is left to the reader.

\begin{lemma}
For a taut contact sphere $(\alpha_1,\alpha_2,\alpha_3)$, the following
conditions are equivalent:
\begin{itemize}
\item[(i)] $\alpha_i\wedge d\alpha_j=0$ for $i\neq j$.
\item[(ii)] The Reeb vector fields $(R_1,R_2,R_3)$ constitute a frame
dual to the coframe $(\alpha_1,\alpha_2,\alpha_3)$.
\item[(iii)] $\beta=0$.
\item[(iv)] The short metric $g_s$ equals the long metric~$g_l$.
\end{itemize}

If any of these conditions holds and $c$ is the normalisation constant,
then one has $[R_i,R_j]=-cR_k$. This implies that the metric $\alpha_1^2+
\alpha_2^2+\alpha_3^2$ has constant curvature
equal to~$c^2/4$. In particular, the metric $g_s=g_l$ has
curvature identically equal to~$1/4$.\qed
\end{lemma}

In analogy with our
terminology in \cite{gego95} we call a taut contact sphere
satisfying any of the conditions in this lemma a
{\bf Cartan structure}.

\begin{remark}
{\rm
Of the taut contact spheres in the statement of
Theorem~\ref{thm:class} exactly those with $\nu =0$ are
Cartan structures. Those with $\nu\neq 0$ are not even
conformally equivalent to a Cartan structure; this follows
for instance from~\cite[Prop.~6.1]{gego95}.
}
\end{remark}

We now want to relate such
Cartan structures to $3$-Sasakian structures. Recall the
definition of these structures.

\begin{definition}
A metric $\og$ on a $3$-manifold $U$ is called {\bf $3$-Sasakian}
if the cone metric $C=e^{2s}(\og+ds^2)$ on $U\times\R$ is
a hyperk\"ahler metric.
\end{definition}

Given such a cone metric $C$, we observe that $\partial_s$ is
a homothetic vector field and $\partial_s\ip C=e^{2s}\, ds$
is a closed $1$-form.

\begin{remark}
{\rm
Given a Riemannian metric $g$ and a vector field $Z$,
the following are equivalent:
\begin{itemize}
\item[(i)] $\nabla Z$ is the identity on every tangent space.
\item[(ii)] $L_Zg=g$ and $Z\ip g$ is closed.
\item[(iii)] $L_Zg=g$, and where $Z$ is non-vanishing it is
orthogonal to a codimension~$1$ foliation.
\item[(iv)] On the open set $\{ Z\neq 0\}$
we have local descriptions $g=e^{2s}(\og+ds^2)$ with $\partial_s=Z$
and $\og$ being the metric induced by $g$ in a transversal orthogonal
to~$Z$.
\end{itemize}
}
\end{remark}

Following~\cite{giry98}, we call a vector field satisfying
either of these conditions a {\bf dilatation}. Thus, any
particular description of $g$ as a cone metric corresponds to
a non-vanishing dilatation.

\begin{proposition}
A $1$-normalised taut contact sphere on a $3$-dimensional manifold $U$
is a Cartan structure if and only if
the induced hyperk\"ahler metric $g$  on $U\times\R$ is a cone metric
with the induced tri-Liouville vector field~$\partial_t$
as dilatation.
\end{proposition}

\begin{proof}
By formula~(\ref{eqn:metric-normal}) for the metric $g$,
the $1$-form $\partial_t\ip g$ equals $e^t(dt+\beta )$,
which is closed if and only if $\beta\equiv 0$.
\end{proof}

\begin{lemma}
Given any metric $\og$ on~$U$,
any endomorphism field on $U\times\R$ parallel with respect to
the corresponding cone metric is invariant under the flow of~$\partial_s$.
\end{lemma}

\begin{proof}
Let $(x_1,x_2,x_3)$ be local coordinates on~$U$. Computing the Levi-Civita
connection of the cone metric in the coordinates $(x_1,x_2,x_3,s)$,
one finds that $e^{-s}\partial_{x_1}, e^{-s}\partial_{x_2},
e^{-s}\partial_{x_3}, e^{-s}\partial_s$ are parallel along the radii.
The coefficients of a parallel endomorphism field in this
frame are constant along the radii. Those coefficients are the same
in the frame $\partial_{x_1},\partial_{x_2},\partial_{x_3},\partial_s$.
\end{proof}

It is a fact in Riemannian geometry that a family of metrics
of the form $\og_{\lambda}=\lambda\og$, $\lambda\in\R^+$, in general gives 
rise to a family of non-isometric cone metrics.
We use the methods of this paper
to give a simple proof of the following result, well-known in Sasakian
geometry, which can be read as saying that for $\og$
having constant positive curvature, only the $\og_{\lambda}$
of curvature equal to $1$ gives rise to a hyperk\"ahler cone.

\begin{proposition}
A $3$-Sasakian metric in dimension $3$ has constant curvature equal to~$1$.
Therefore, if a $4$-dimensional hyperk\"ahler metric admits a dilatation,
it must be flat.
\end{proposition}

\begin{proof}
Let $\og$ be a Riemannian metric on a $3$-manifold $U$ and assume
that the cone metric $C=e^{2s}(ds^2+\og)$ is hyperk\"ahler.
By the preceding lemma the
flow of the dilatation $(1/2)\partial_s$ is tri-holomorphic and also
tri-Liouville. Introduce the new coordinate $t=2s$, so that
$\partial_t=(1/2)\partial_s$. By the theory we have
developed, there is a $1$-normalised taut contact sphere
$(\alpha_1,\alpha_2,\alpha_3)$  on $U$ such that
\[ C=e^t(dt^2+\alpha_1^2+\alpha_2^2+\alpha_3^2).\]
Notice that $\og$ is recovered from the cone metric and
the dilatation $\partial_s$ via the formula
$\og =C|_{\{ s=0\}}=C|_{\{ \|\partial_s\| =1\}}$. But
$\|\partial_s\| =2\|\partial_t\|$, so
\[ \og = e^t(dt^2+\alpha_1^2+\alpha_2^2+\alpha_3^2)|_{\{e^t=1/4\}}
=(\alpha_1/2)^2+(\alpha_2/2)^2+(\alpha_3/2)^2.\]
We see that $\og$ has a $2$-normalised Cartan structure as an
orthonormal frame, making it a metric of constant curvature
equal to~$1$.

Then $e^{2s}(\og + ds^2)$ describes the
$4$-dimensional Euclidean metric in spherical coordinates.
\end{proof}

This result must be read as {\em local} rigidity of $3$-Sasakian structures.
As our constructions of taut contact spheres giving rise to non-flat
hyperk\"ahler metrics show, taut contact spheres have much
richer local geometry. In particular, their associated hyperk\"ahler
metrics admit homotheties, but no dilatations. So they
are not cone metrics in any way.

As we have mentioned in Section~\ref{section:class}, the
$c$-normalised taut contact spheres on closed $3$-manifolds are
orthonormal for Berger metrics, which are spherical only in the
Cartan case. The long metrics, on the other hand, are always
spherical.

\vspace{2mm}

By the discussion in this section,
our Theorems \ref{thm:classold}
and~\ref{thm:class} can be read as a
classification of the closed $3$-Sasakian $3$-manifolds:

\begin{corollary}
\label{cor:Sasaki}
The closed
$3$-Sasakian $3$-manifolds are precisely the left-quotients
of~$\SU (2)$. \qed
\end{corollary}

For the parametric family of taut contact spheres
in Theorem~\ref{thm:class}, the induced hyperk\"ahler metric
is always the standard Euclidean metric on $\E^4\setminus\{ 0\}$,
so it is always a cone metric. But only for the parameter
value $\nu =0$ (corresponding to the unique class containing Cartan
structures) does the tri-Liouville vector field $Y$ point in
the radial direction (i.e.\ the direction of the only non-vanishing
dilatation of $\E^4\setminus\{ 0\}$).

In fact, only the tri-Liouville vector field changes within this
parametric family.
The $d\alpha_i$ obviously do not depend on~$\nu$,
and a little computation shows that neither do the Reeb
vector fields~$R_i$ of the~$\alpha_i$.

A proof of Corollary~\ref{cor:Sasaki} previous to the one
given above was indeed based on the observation that a
$3$-Sasakian $3$-manifold is a space of constant curvature~$1$;
the classification of $3$-Sasakian manifolds among
the $3$-dimensional space forms was achieved by
Sasaki~\cite{sasa72} (who still spoke of {\em normal contact
metric $3$-structures}\/).

\section*{Appendix: A contact transformation leading to Helmholtz' equation}
The following intriguing proposition relates solutions of
equation~(\ref{eqn:MA2}) from Section~\ref{section:Helmholtz}
to solutions of the Helmholtz equation (cf.~\cite{evan98} for this
terminology) $\Delta u+2u=0$ on the $2$-sphere $S^2$ with its
metric of constant curvature~$1$. In geodesic polar coordinates this metric
is given by
\[ ds^2=dr^2+\sin^2r\, d\theta^2.\]
Hence the gradient of a differentiable
function $f\co S^2\rightarrow \R$
is computed via
\[ \nabla f=f_r\partial_r+\frac{1}{\sin^2r}f_{\theta}\partial_{\theta}.\]
The area element is $A=\sin r\, dr\wedge d\theta$, so from
\[ d(X\ip A)=L_XA=\mbox{\rm div} (X)\cdot A\;\;\mbox{\rm and}\;\;
\Delta f=\mbox{\rm div} (\nabla f)\]
we see that the {\it spherical} Laplacian $\Delta =\Delta^{S^2}$ takes
the form
\[ \Delta f=\frac{\cos r}{\sin r}f_r+f_{rr}+
\frac{1}{\sin^2r}f_{\theta\theta}.\]

\begin{proposition}
\label{prop:helmholtz}
Let $h(s_1,s_2)$ be a solution of equation~$(\ref{eqn:MA2})$. Set
\[ (T)\;\left\{ \begin{array}{rcl}
\theta & = & s_1,\\
r & = & \mbox{\rm arccot} (h_{s_2})\in (0,\pi ),\\
u & = & -s_2\cos r+h \sin r.
\end{array}\right. \]
If $h_{s_2s_2}\neq 0$, then $u=u(r,\theta )$ is a function of the
independent variables $r$ and $\theta$ that solves the
spherical Helmholtz equation
\begin{equation}
\label{eqn:Helmholtz}
\Delta u+2u=0.
\end{equation}

Conversely, a solution $u(r,\theta )$ of $(\ref{eqn:Helmholtz})$
satisfying $u+u_{rr}\neq 0$
gives rise to independent variables $s_1,s_2$ and a solution
$h(s_1,s_2)$ of $(\ref{eqn:MA2})$ via the inverse transformation
\[ (T^{-1})\;\left\{ \begin{array}{rcl}
s_1 & = & \theta,\\
s_2 & = & u_r\sin r-u\cos r,\\
h   & = & u_r\cos r+u\sin r.
\end{array}\right. \]
\end{proposition}

\begin{proof}
Given $h(s_1,s_2)$, let $(r,\theta ,u)$ be defined by~$(T)$. The condition
$h_{s_2s_2}\neq 0$ is obviously necessary and sufficient for $(s_1,s_2)\mapsto
(r,\theta )$ to be an invertible coordinate transformation. From
\begin{eqnarray*}
du & = & -\cos r\, ds_2+s_2\sin r\, dr+(h_{s_1}\, ds_1+h_{s_2}\, ds_2)\sin r
+h\cos r\, dr\\
  & = & (s_2\sin r+h\cos r)\, dr+h_{s_1}\sin r d\theta
\end{eqnarray*}
we see that $u$ is a function of $r$ and $\theta$ with
\[ (T)'\;\left\{\begin{array}{rcl}
u_r & = & s_2\sin r+h\cos r,\\
u_{\theta} & = & h_{s_1}\sin r.
\end{array}\right. \]

Conversely, given $u(r,\theta )$, let $(s_1,s_2,h)$ be defined by~$(T^{-1})$.
Since
\[ ds_2=(u+u_{rr})\sin r\, dr+(u_{r\theta}\sin r-u_{\theta}\cos r)\,
d\theta ,\]
the condition $u+u_{rr}\neq 0$ is necessary and sufficient for
$(r,\theta )\mapsto (s_1,s_2)$ to be a coordinate transformation,
because $r\in (0,\pi )$. From
\begin{eqnarray*}
dh & = & (u+u_{rr})\cos r\, dr+(u_{r\theta}\cos r+u_{\theta}\sin r)\, d\theta\\
  & = & \frac{u_{\theta}}{\sin r}\, ds_1+\cot r\, ds_2
\end{eqnarray*}
we infer that $h$ is a function of $s_1$ and $s_2$ with
\[ (T^{-1})'\;\left\{\begin{array}{rcl}
h_{s_1} & = & u_{\theta}/\sin r,\\
h_{s_2} & = & \cot r.
\end{array}\right. \]

It is then a straightforward check that $(T)$ and $(T^{-1})$ are indeed
inverse transformations of each other. In particular, one finds that
\[ h_{s_2s_2}=-\frac{1}{(u+u_{rr})\sin^3r},\]
which shows that $(T)$ is defined on $h$ if and only if $(T^{-1})$
is defined on~$u$.

In order to show that solutions of (\ref{eqn:MA2}) correspond
to solutions of~(\ref{eqn:Helmholtz}) under the transformation~$(T)$,
it is convenient to rewrite
(\ref{eqn:MA2}) as an exterior differential
system
\begin{equation}
\label{eqn:system}
\begin{array}{rcl}
dh & = & p_1\, ds_1+p_2\, ds_2,\\
ds_1\wedge ds_2 & = & dp_1\wedge dp_2+h\, ds_1\wedge dp_2-
                      p_2^2\, ds_1\wedge ds_2.
\end{array}
\end{equation}
We then compute, using $(T^{-1})$ and $(T^{-1})'$,
\begin{eqnarray*}
\lefteqn{dp_1\wedge dp_2+h\, ds_1\wedge dp_2-(1+p_2^2)\, ds_1\wedge ds_2=}\\
  & = & d\frac{u_{\theta}}{\sin r}\wedge d\cot r+(u_r\cos r+u\sin r)\,
        d\theta\wedge d\cot r\\
  &   & \mbox{}-(1+\cot^2r)\, d\theta\wedge d(u_r\sin r-u\cos r)\\
  & = & \left( \frac{u_{\theta\theta}}{\sin r}\cdot \frac{1}{\sin^2r}
        + u_r\, \frac{\cos r}{\sin^2 r}+u\, \frac{1}{\sin r}
        +\frac{1}{\sin^2 r}\, (u+u_{rr})\sin r\right)\,
         dr\wedge d\theta\\
  & = & \frac{1}{\sin r}\left( \frac{1}{\sin^2r}\, u_{\theta\theta}+
        \frac{\cos r}{\sin r}\, u_r+u_{rr}+2u\right)\, dr\wedge d\theta\\
  & = & \frac{1}{\sin r}\, (\Delta u+2u)\, dr\wedge d\theta. 
\end{eqnarray*}
This completes the proof of Proposition~\ref{prop:helmholtz}.
\end{proof}

Since the reader is bound to wonder how we arrived at the transformation~$(T)$,
we present, {\em in nuce}, our chain of discovery:
First, one can get rid of mixed derivatives in the exterior differential
system~$(\ref{eqn:system})$ by introducing
$p_2$ as an independent variable.
The identity
\[ dh-p_1\, ds_1-p_2\, ds_2=d(h-p_2s_2)-p_1\, ds_1+s_2\, dp_2\]
suggests the contact transformation $h(s_1,s_2)\rightsquigarrow k(t_1,t_2)$
given by
\[ t_1=s_1,\; t_2=p_2=h_{s_2},\; k=h-s_2h_{s_2},\]
with inverse transformation
\[ s_1=t_1,\; s_2=-k_{t_2},\; h=k-t_2k_{t_2}.\]
(Here one needs $h_{s_2s_2}\neq 0$ or $k_{t_2t_2}\neq 0$, respectively,
for these to be honest, i.e.\ invertible, transformations.)
We compute
\begin{eqnarray*}
dp_1\wedge dp_2 & = & dk_{t_1}\wedge dt_2 \; = \; k_{t_1t_1}\, dt_1\wedge dt_2,\\
h\, ds_1\wedge dp_2 & = & (k-t_2k_{t_2})\, dt_1\wedge dt_2,\\
-(1+p_2^2)\, ds_1\wedge ds_2 & = & (1+t_2^2)\, dt_1\wedge dk_{t_2}
                                   \; = \; (1+t_2^2)k_{t_2t_2}\,
                                            dt_1\wedge dt_2.
\end{eqnarray*}
So the contact transformation takes equation (\ref{eqn:MA2})
to the following equation,
where we now write $k_i$ for $k_{t_i}$ etc.:
\begin{equation}
\label{eqn:k}
k_{11}+(1+t_2^2)k_{22}-t_2k_2+k=0.
\end{equation}

The first order term in this equation can be made to disappear by
setting $x=t_1$, $\sinh y=t_2$ and $w(x,y)=k(x,\sinh y)/\cosh y$. This
turns equation $(\ref{eqn:k})$ into
\begin{equation}
\label{eqn:w}
w_{xx}+w_{yy}+\frac{2}{\cosh^2y}w=0,
\end{equation}
which is the Helmholtz equation for the metric
\[ \left(\begin{array}{cc}
\frac{1}{\cosh^2y} & \; 0\\
0                  & \; \frac{1}{\cosh^2y}
\end{array}\right) . \]
The Gau{\ss} curvature of this metric turns out to be identically
equal to~$1$. Indeed, it is the spherical metric in Mercator
coordinates. Therefore, we pass from (\ref{eqn:w}) to
(\ref{eqn:Helmholtz}) by making
the substitution
\[ \cosh y =1/\sin r,\;\; \sinh y =\cot r,\;\;
\theta =x,\;\;\mbox{\rm and}\;\; u(r,\theta )=w(x,y).\]

\vspace{2mm}

We now want to find an explicit solution of~$(\ref{eqn:Helmholtz})$.
The ansatz
\[ u(r,\theta )=\cos\theta \sin r\, g(r)\]
leads to
\[ 3g_r\cos r+g_{rr}\sin r =0,\]
which has the solution $g_r=1/\sin^3r$; that in turn integrates to
\[ g(r)=\int_{\cot r}^0\sqrt{1+t^2}\, dt.\]
The resulting $u$ corresponds under $(T^{-1})$ to the solution
\[ h(s_1,s_2)=\cos s_1\, \int_{1}^{s_2/\cos s_1} \sqrt{\xi^2-1}\, d\xi\]
of equation~$(\ref{eqn:MA2})$.
As domain of definition we may take
\[ U'=\{ (s_1,s_2)\in \R^2\co |s_1|<\pi /2,\; s_2>\cos s_1\} .\]

The K\"ahler potential $H$ corresponding to this solution
gives rise to a hyperk\"ahler metric~$g$ on $U\times \R$
inducing a taut contact sphere on~$U$, with
\[ U=\{ (s_1,s_2,s_3)\in \R^3\co (s_1,s_2)\in U'\}.\]
The relation with the complex coordinates is given by $t+\bi
s_1=2z_1$ and $s_2+\bi s_3=2z_2$, say.

\vspace{2mm}

We claim that $g$ is non-flat.
The coefficients of this metric are
\[ g_{\alpha\obeta}=g(\partial_{z_{\alpha}},\partial_{\oz_{\beta}})=
  \Omega_1(\partial_{z_{\alpha}},J_1\partial_{\oz_{\beta}})=
  \frac{1}{2}H_{z_{\alpha}\oz_{\beta}}.\]
Write $\bG$ for the $(2\times 2)$-matrix $(g_{\alpha\obeta})$. Then the
curvature tensor $K_{\alpha\obeta\gamma\odelta}$, read for fixed
$\gamma ,\odelta$ as a $(2\times 2)$-matrix indexed by $\alpha$
and $\obeta$, is computed by
\[ (K_{\alpha\obeta\gamma\odelta})=\bG_{z_{\gamma}\oz_{\delta}}
-\bG_{z_{\gamma}}\bG^{-1}\bG_{\oz_{\delta}},\]
cf.~\cite[p.~159]{kono69}.

We now want to show that $K_{2\otwo 2\otwo}$ is non-zero. In the following
computations we write $\ast$ for any matrix entry that is irrelevant for the
final result:
\begin{eqnarray*}
\bG & = & \frac{1}{2}e^{z_1+\oz_1} \left( \begin{array}{cc}
             h+h_{11}       & \; h_2-\bi h_{12}\\
             h_2+\bi h_{12} & \; h_{22}
            \end{array}\right) ,\\
\bG_{z_2} & = & \frac{1}{2}e^{z_1+\oz_1} \left( \begin{array}{cc}
             \ast               & \; \ast\\
             h_{22}+\bi h_{122} & \; h_{222}
            \end{array}\right) ,\\
\bG_{\oz_2} & = & \frac{1}{2}e^{z_1+\oz_1} \left( \begin{array}{cc}
             \ast & \; h_{22}-\bi h_{122}\\
             \ast & \; h_{222}
            \end{array}\right) ,\\
\bG_{z_2\oz_2} & = & \frac{1}{2}e^{z_1+\oz_1} \left( \begin{array}{cc}
             \ast & \; \ast\\
             \ast & \; h_{2222}
            \end{array}\right) .
\end{eqnarray*}
On the hyperplane $\{ 2\,\mbox{\rm Im}(z_1)=z_1-\oz_1=0\}$, corresponding to
the line $\{ s_1=0\}$, we have
\[ h_1=0 \;\;\mbox{\rm and}\;\; h=\int_{1}^{s_2}\sqrt{\xi^2-1}\, d\xi .\]
Writing $\sqrt{s_2^2-1}=\sigma$, for short, we have along that same line
\[ h_2=\sigma ,\; h_{22}=s_2/\sigma ,\; h_{222}=-1/\sigma^3,\;
  h_{2222}=3s_2/\sigma^5.\]
At $(s_1,s_2)=(0,\sqrt{2})$ we thus find
\[ h_2=1,\; h_{22}=\sqrt{2},\; h_{222}=-1,\; h_{2222}=3\sqrt{2},\]
\[ h_{12}=h_{122}=0,\; \mbox{\rm and}\; h+h_{11}=\sqrt{2}.\]
That last equality can be computed from (\ref{eqn:MA2}). Hence, at $(z_1,z_2)=
(0,\sqrt{2}/2)$ we obtain
\[ \left( \begin{array}{cc}
   \ast & \; \ast\\
   \ast & \; K_{2\otwo 2\otwo}
   \end{array} \right) =
\frac{1}{2}\left( \begin{array}{cc}
   \ast & \; \ast\\
   \ast & \; 3\sqrt{2}
   \end{array} \right)
-\frac{1}{2}\left( \begin{array}{cc}
   \ast     & \; \ast\\
   \sqrt{2} & \; -1
   \end{array} \right)
\left( \begin{array}{cc}
   \sqrt{2} & \; -1\\
   -1       & \; \sqrt{2}
   \end{array} \right)
\left( \begin{array}{cc}
   \ast & \; \sqrt{2}\\
   \ast & \; -1
   \end{array} \right) ,\]
which yields $K_{2\otwo 2\otwo}=-\sqrt{2}\neq 0$. Thus $g$ is indeed
non-flat.
\begin{acknowledgements}
We are immensely grateful to Nigel Hitchin for drawing our attention
to the similarities between our original treatment of taut
contact spheres and the Gibbons-Hawking ansatz. Hitchin's observation
inspired a large part of the present paper.
We thank Martin L\"ubke for useful hints concerning the curvature
calculations in the appendix.
\end{acknowledgements}




\begin{thebibliography}{99}
%
\bibitem{ahs78}
Atiyah, M. F.,  Hitchin, N. J., and Singer, I. M.:
Self-duality in four-dimensional Riemannian geometry,
{\it Proc. Roy. Soc. London Ser. A}
{\bf 362}, 425--461 (1978)
%
\bibitem{basf97}
Bakas, I. and Sfetsos, K.: Toda fields of $\SO (3)$ hyper-K\"ahler
metrics and free field realizations, {\it Internat. J. Modern
Phys. A} {\bf 12}, 2585--2611 (1997)
%
\bibitem{bar93}
B\"ar, C.: Real Killing spinors and holonomy, {\it Comm. Math.
Phys.} {\bf 154}, 509--521 (1993)
%
\bibitem{bess87}
Besse, A.:
{\it Einstein Manifolds},
Ergeb. Math. Grenzgeb. (3).
Berlin: Springer-Verlag, 1987
%
\bibitem{boye88}
Boyer, C. P.:
A note on hyper-Hermitian four-manifolds,
{\it Proc. Amer. Math. Soc.} {\bf 102}, 157--164 (1988)
%
\bibitem{boga99}
Boyer, C. P. and Galicki, K.:
$3$-Sasakian manifolds.
In: {\it Surveys in Differential Geometry: Essays on Einstein Manifolds},
Surv.  Differ. Geom. {\bf VI}.
Boston: Int. Press, 1999, pp. 123--184
%
\bibitem{boga08}
Boyer, C. P. and Galicki, K.:
{\it Sasakian Geometry}.
Oxford Math. Monogr.
Oxford: Oxford University Press, 2008
%
\bibitem{cala58}
Calabi, E.: Improper affine hyperspheres of convex type and a
generalization of a theorem by K.~J\"orgens,
{\it Michigan Math.  J.} {\bf 5}, 105--126 (1958)
%
\bibitem{cala70}
Calabi, E.: Examples of Bernstein problems for some nonlinear
equations.
In: {\it Global Analysis} (Berkeley, 1968),
Proc. Sympos. Pure Math. {\bf 15}.
Providence: American Mathematical Society, 1970, pp. 223--230
%
\bibitem{cala75}
Calabi, E.: A construction of nonhomogeneous Einstein metrics,
in: {\it Differential Geometry} (Stanford, 1973),
Proc. Sympos. Pure Math. {\bf 27}, Part~2.
Providence: American Mathematical Society, 1975, pp. 17--24
%
\bibitem{ctv96}
Chave, T., Tod, K. P., and Valent, G.:
$(4,0)$ and $(4,4)$ sigma models with a tri-holomorphic
Killing vector,
{\it Phys. Lett. B} {\bf 383}, 262--270 (1996)
%
\bibitem{evan98}
Evans, L. C.:
{\it Partial Differential Equations},
Grad. Stud. Math. {\bf 19}.
Providence: American Mathematical Society, 1998
%
\bibitem{geig96}
Geiges, H.:
Symplectic couples on $4$-manifolds,
{\it Duke Math. J.} {\bf 85}, 701--711 (1996)
%
\bibitem{geig97}
Geiges, H.:
Normal contact structures on $3$-manifolds,
{\it T\^{o}hoku Math. J.} {\bf 49}, 415--422 (1997)
%
\bibitem{gego95}
Geiges, H. and Gonzalo, J.:
Contact geometry and complex surfaces,
{\it Invent. Math.} {\bf 121}, 147--209 (1995)
%
\bibitem{gego97}
Geiges, H. and Gonzalo, J.:
Contact circles on $3$-manifolds,
{\it J. Differential Geom.} {\bf 46}, 236--286 (1997)
%
\bibitem{gego02}
Geiges, H. and Gonzalo, J.:
Moduli of contact circles,
{\it J. Reine Angew. Math.}
{\bf 551}, 41--85 (2002)
%
\bibitem{giha78}
Gibbons, G. W. and Hawking, S. W.:
Gravitational multi-instantons,
{\it Phys. Lett. B} {\bf 78}, 430--432 (1978)
%
\bibitem{giry98}
Gibbons, G. W. and Rychenkova, P.:
Cones, tri-Sasakian structures and superconformal invariance,
{\it  Phys. Lett. B} {\bf 443}, 138--142 (1998)
%
\bibitem{gkn00}
Godli\'nski, M., Kopczy\'nski, W., and Nurowski, P.:
Locally Sasakian mani\-folds,
{\it Classical Quantum Gravity}
{\bf 17}, L105--L115 (2000)
%
\bibitem{kono69}
Kobayashi, S. and Nomizu, K.:
{\it Foundations of Differential Geometry~II}.
New York: Interscience, 1969
%
\bibitem{sasa72}
Sasaki, S.:
Spherical space forms with normal contact metric
$3$-structure,
{\it J. Differential Geom.}
{\bf 6}, 307--315 (1972)
%
%
%
\end{thebibliography}
\end{document}